\newif\ifdraft
\draftfalse

\documentclass[reqno,a4paper]{amsart}

\usepackage{mathrsfs,amsmath,amssymb,graphicx}
\usepackage{stmaryrd}
\usepackage[unicode]{hyperref}
\usepackage{latexsym}
\usepackage{layout}
\usepackage[english]{babel}
\usepackage{color}
\usepackage{subcaption}
\usepackage{fancyvrb}
\usepackage{color}
\usepackage{amsmath,amscd}
\usepackage{multicol}
\usepackage{multirow}
\usepackage{enumitem}
\usepackage{overpic}

\usepackage{mathtools}
\usepackage{tikz} 
 
\usetikzlibrary{matrix,arrows,decorations,shapes,calc}
\tikzset{every picture/.append style={remember picture},
na/.style={baseline=-.5ex}}

\newif\ifhknots
\newif\ifhlinks
\newif\ifgenusthree
\hknotsfalse
\hlinksfalse
\genusthreetrue

\theoremstyle{plain}
\newtheorem{theorem}{Theorem}[section]
\newtheorem{lemma}[theorem]{Lemma}

\newtheorem{corollary}[theorem]{Corollary}

\theoremstyle{definition}
\newtheorem{definition}{Definition}[section]
\theoremstyle{remark}
\newtheorem{remark}{Remark}[section]

\newcommand{\nada}[1]   {}

\newcommand{\Dcal} {\mathcal{D}}
\newcommand{\Done} {\Dcal_1}
\newcommand{\Dtwo} {\Dcal_2}
\newcommand{\DtwoA} {\Dcal_2^a}
\newcommand{\DtwoB} {\Dcal_2^b}
\newcommand{\Dthree} {\Dcal_3}

\newcommand{\numcc}{C}

\newcommand{\Stwo} {\mathfrak S}

\newcommand{\HL}{{\rm HL}}

\newcommand{\Sbb} {\mathbb S}

\newcommand{\sphere} {\Sbb^3}

\definecolor{mygray}{rgb}{0.92,0.92,0.92}

\newcommand{\Compl}[1]{E(#1)}%W
\newcommand{\fthree}{\op{\mathfrak{f}3}}
\newcommand{\ftwo}{\op{\mathfrak{f}2}}
\newcommand{\flink}{\op{\ell}}

\newcommand{\onesum}{\operatorname{-{}-}}
\newcommand{\twosum}{\operatorname{\#}}
\newcommand{\op}[1]{\operatorname{#1}}

\newcommand{\knotlink}{knot/link}
\newcommand{\knotlinks}{knot/links}

\newcommand{\draftMP}[1]{\ifdraft{\color{blue}#1}\fi}
\newcommand{\draftGB}[1]{\ifdraft{\color{red}#1}\fi}
\newcommand{\draftYW}[1]{\ifdraft{\color{orange}#1}\fi}
\newcommand{\draftGP}[1]{\ifdraft{\color{teal}#1}\fi}

\graphicspath{{figures/}}

\numberwithin{equation}{section}
\numberwithin{figure}{section}

\title[Handlebody-links with Euler characteristic $-2$]{Handlebody-links with Euler characteristic $-2$ up to six crossings
}
\author[G.\ Bellettini]{Giovanni Bellettini}
\address{Dipartimento di Scienze Matematiche, Informatiche e Fisiche,
Universit\`a di Udine, Italy,
e-mail: giovanni.bellettini@uniud.it
and International Centre for Theoretical Physics ICTP,
Mathematics Section, 34151 Trieste, Italy
}
\email{giovanni.bellettini@unisi.it}

\author[G.\ Paolini]{Giovanni Paolini}
\address{Dipartimento di Matematica, Universit\`a di  Bologna, Italy}
\email{g.paolini@unibo.it}

\author[M.\ Paolini]{Maurizio Paolini}
\address{Dipartimento di Matematica e Fisica, Universit\`a Cattolica del Sacro Cuore, Brescia, Italy}
\email{maurizio.paolini@unicatt.it}

\author[Y.\ S.\ Wang]{Yi-Sheng Wang}
\address{National Sun Yat-sen University, Kaohsiung, Taiwan}
\email{yisheng@math.nsysu.edu.tw}

\date{\today}

\begin{document}

%{\tiny \tableofcontents}
\thanks{}

\begin{abstract}
We enumerate all non-split
handlebody-links  
with Euler characteristic $-2$
up to six crossings, up to ambient isotopy and mirror image.
\end{abstract}

\maketitle

%%%%%%%%%%%%%%%%%%%%%%%%%%%%%%%%%%%%%%%%%%%%%%%%%%%%%%%%%%%%%%%%%%%%%%%%%%%%%
%
% for any new section simply add the corresponding two lines here
% and modify accordingly the Makefile
%
%%%%%%%%%%%%%%%%%%%%%%%%%%%%%%%%%%%%%%%%%%%%%%%%%%%%%%%%%%%%%%%%%%%%%%%%%%%%%

\draftMP{DraftMP comments by Maurizio}
\draftYW{DraftYW comments by Yi-Sheng}
\draftGB{DraftGB comments by Giovanni B}
\draftGP{DraftGP comments by Giove}

%%%%%%%%%%%%%%%%%%%%%%%%%%%%%%%%%%%%%%%%%%%%%%%%%%%%%%%%%%%%%%%%%%%%%%%%%%%%%%
\section{Introduction}\label{sec:intro}

A %(unoriented) 
\emph{spatial graph} $\Gamma$ is a finite, not necessarily connected, graph embedded in $S^3$. In particular, classical links are spatial graphs with Euler characteristic $\chi=0$, and a classical knot is a one-component classical link. 
A spatial graph is 
\emph{trivalent} if all its vertices are 
trivalent, and a trivalent spatial graph with Euler characteristic $\chi$ has $-2\chi$ trivalent vertices. 

A \emph{handlebody-link} $V$ is a finite union of 
disjoint handlebodies of positive genus embedded in $S^3$. A \emph{handlebody-knot} $V$ is a one-component handlebody-link, and the genus of the component is the genus of $V$; in this case, we have $g=1-\chi$. A \emph{multi-component} handlebody-link is a handlebody-link with more than one component.
%In general, we say 
%
%\begin{equation}\label{eq:h_type}
%V {\rm ~ has~ htype~}
%[n_1,\cdots, n_k]
%\end{equation}
%
%if the number of 
%genus $i\geq 1$ components in $V$ is $n_i$. 
%In particular, $V$ has $\sum_{i=1}^k n_i$ components, and the Euler characteristic $\chi$ of $V$ is given by
%(having $\sum_{i=1}^k n_i$ 
%connected components)
%can be computed by 
%$\chi=\sum_{i=1}^k (1-i)n_i$. 

A handlebody-link $V$ is \emph{non-split} if every $2$-sphere $S$ with $S\cap V=\emptyset$ bounds a $3$-ball disjoint from $V$. A non-split handlebody-link $V$ is \emph{reducible} if there exists a $2$-sphere $S$ meeting $V$ in an essential disk, and is \emph{irreducible} otherwise.

Two handlebody-links
%\knotlinks\ 
$V,V'$ (resp.\ spatial graphs $\Gamma,\Gamma'$) are \emph{equivalent} if there is an orientation-preserving self-homeomorphism $f$ of $S^3$ such that $f(V)=V'$ (resp.\ $f(\Gamma)=\Gamma'$). Two handlebody-links $V,V'$ are equivalent, \emph{up to mirror image}, 
if $V$ is equivalent to $V'$ or the mirror image of $V'$.
While a regular neighborhood of a spatial graph determines a handlebody-link, up to equivalence, a handlebody-link  may have infinitely many spines, inequivalent as spatial graphs. On the other hand, a solid torus has a unique spine, up to ambient isotopy, so the study of handlebody-link
with $\chi=0$ is equivalent to classical link theory. For $\chi<0$, every handlebody-link
admits a trivalent spine, and it is a standard convention to represent a handlebody-link by one of its trivalent spines.

A \emph{diagram} of a handlebody-link
%\knotlink\ 
is a diagram of one of its spines, and the \emph{crossing number} of a handlebody-link
%\knotlink\ 
is the minimal number of crossings among all its diagrams.
%Genus two handlebody-knots, or handlebody-knots with $\chi=-1$, up to six crossings are enumerated in Ishii-Kishimoto-Moriuchi-Suzuki 
The tables in Ishii-Kishimoto-Moriuchi-Suzuki
\cite{IshKisMorSuz:12} and Bellettini-Paolini-Paolini-Wang \cite{BePaPaWa:23} together enumerate 
all irreducible non-split handlebody-links with $\chi=-1$, 
up to six crossings, up to equivalence and mirror image.
%while $n$-component handlebody-links with $n>1$ and $\chi=-1$
%of htype
% $[n,1]$, $n >0$,
% up to six crossings, are enumerated in \cite{BePaPaWa:23}. They together enumerate all handlebody-\knotlinks\ of $\chi=-1$ up to six crossings. 
%In \cite{BePaPaWa:25} all genus two handlebody-knots with seven crossings have been
%enumerated (up to mirror image),
%thus extending the 
%table in \cite{IshKisMorSuz:12}.
%
The main results here extend their works to the case of $\chi=-2$.  
\begin{theorem}
Table \ref{tab:hknots} enumerates, up to six crossings, all \textit{irreducible} genus $3$ handlebody-knots, namely handlebody-knots with $\chi=-2$, up to equivalence and mirror image.  
\end{theorem}
This complements the genus two handlebody-knot table in \cite{IshKisMorSuz:12}.
%A handlebody-\knotlink\ is irreducible if it is not reducible in the sense of Definition \ref{def:reducible}.

\begin{theorem}
Table \ref{tab:hlinks} enumerates, up to six crossings, all irreducible non-split multi-component handlebody-links with $\chi=-2$, up to equivalence and mirror image.  
\end{theorem}
This complements the handlebody-link table in \cite{BePaPaWa:23}. 
%A handlebody-link with $\chi=-2$ may contain two genus two components or one genus three component, together with some solid tori. In general, 
To encode the genus of each component in a handlebody-link $V$, we say  
\begin{equation}\label{eq:h_type}
V {\rm ~ is~ of~ type~}
[n_1,\cdots, n_k]
\end{equation}
if $V$ has $n_i$ components of 
genus $i\geq 1$. 
Since the Euler characteristic $\chi$ of $V$ is given by 
$\chi=\sum_{i=1}^k (1-i)n_i$, a 
handlebody-link with $\chi=-1$ is of type $[n,1]$, and a 
handlebody-link with $\chi=-2$ is either of type $[n,0,1]$, 
or $[n,2]$, $n\in \mathbb{N}$. 
In Table \ref{tab:hlinks}, handlebody-links $4_1$ to $6_4$ are 
of type $[0,2]$, 
$6_5$ to $6_7$ are of type $[1,0,1]$, 
$6_8$ to $6_{10}$ are of type $[1,2]$ and lastly $6_{11}$ is of type $[2,0,1]$.

In addition, reducible non-split handlebody-links are enumerated in Table \ref{tab:typeone} and Figure \ref{fig:typeonedetail} by their type, up to equivalence and mirror image.
Thus all non-split handlebody-links with $\chi=-2$ are classified.

Section \ref{sec:diagram} recalls some preliminaries on diagrams of handlebody-links and spatial graphs. Section \ref{sec:completeness} proves the completeness of Tables \ref{tab:hknots}, \ref{tab:hlinks}. Sections \ref{sec:uniqueness} and \ref{sec:hardpairs} 
show there is no duplication in the tables, and Section \ref{sec:irreducibility} proves the non-splittability and irreducibility of entries therein.

%%%%%%%%%%%%%%%%%%%%%%%%%%%%%%%%%%%%%%%%%%%%%%%%%%%%%%%%%%%%%%%%%%%%%%%%%%%%%%
\newdimen\entrysize \entrysize=59pt
\ifdraft
\newcommand{\fffbox}{\framebox}
\else
\newcommand{\fffbox}{\makebox}
\fi
\ifdraft
\newcommand{\tabentry}[2]{\fffbox[\entrysize]{\includegraphics[height=0.15\textwidth]{figures/ht001x#1n#2}}$#1_{#2}$}
\else
\newcommand{\tabentry}[2]{\fffbox[\entrysize]{\includegraphics[height=0.15\textwidth]{figures/ht001x#1n#2}}$#1_{#2}$}
\fi
\newcommand{\tabentrym}[2]{\fffbox[\entrysize]{\includegraphics[height=0.15\textwidth]{figures/tmissing}}$#1_{#2}$}
\ifdraft
\newcommand{\tabentryt}[1]{\fffbox[\entrysize]{\begin{tikzpicture}[scale=0.8,every edge/.style={draw,thick}]
\draw[green,thin] (-1.05,-1.05) rectangle (1.05,1.05);
\input{tikz/l#1.tikz}\end{tikzpicture}}$x_{#1}$}
\else
\newcommand{\tabentryt}[1]{\fffbox[\entrysize]{\begin{tikzpicture}[scale=0.8,every edge/.style={draw,thick}]\input{tikz/#1.tikz}\end{tikzpicture}}$x_{#1}$}
\fi

\begin{table}[ht]
\caption{Irreducible handlebody-knots of genus $3$ up to six crossings.}
\label{tab:hknots}
\tabentry{5}{1}%
\tabentry{6}{1}%
\tabentry{6}{2}%
\tabentry{6}{3}%
\tabentry{6}{4}%
\\
\smallskip
\tabentry{6}{5}%
\tabentry{6}{6}%
\tabentry{6}{7}%
\tabentry{6}{8}%
\tabentry{6}{9}%
\ifdraft\tabentry{6}{10old}\fi%
\\
\smallskip
\ifdraft\tabentry{6}{11old}\fi%
\tabentry{6}{10}%
\tabentry{6}{11}%
\tabentry{6}{12}%
\tabentry{6}{13}%
\end{table}

\newdimen\entrysize \entrysize=59pt
\newcommand{\tabentryla}[3]{\fffbox[\entrysize]{\includegraphics[height=0.15\textwidth]{figures/ht02x#1#2}}$#1_{#3}$}
\newcommand{\tabentrylb}[3]{\fffbox[\entrysize]{\includegraphics[height=0.15\textwidth]{figures/ht101x#1#2}}$#1_{#3}$}
\newcommand{\tabentrylc}[3]{\fffbox[\entrysize]{\includegraphics[height=0.15\textwidth]{figures/ht12x#1#2}}$#1_{#3}$}
\newcommand{\tabentryld}[3]{\fffbox[\entrysize]{\includegraphics[height=0.15\textwidth]{figures/ht201x#1#2}}$#1_{#3}$}
%\newcommand{\tabentrym}[2]{\fffbox[\entrysize]{\includegraphics[height=0.15\textwidth]{figures/tmissing}}$#1_{#2}$}
%\ifdraft
%\newcommand{\tabentryt}[1]{\fffbox[\entrysize]{\begin{tikzpicture}[scale=0.8,every edge/.style={draw,thick}]
%\draw[green,thin] (-1.05,-1.05) rectangle (1.05,1.05);
%\input{tikz/l#1.tikz}\end{tikzpicture}}$x_{#1}$}
%\else
%\newcommand{\tabentryt}[1]{\fffbox[\entrysize]{\begin{tikzpicture}[scale=0.8,every edge/.style={draw,thick}]\input{tikz/#1.tikz}\end{tikzpicture}}$x_{#1}$}
%\fi

\begin{table}[ht]
\caption{Irreducible multi-component handlebody-links with Euler characteristic $-2$, up to six crossings.
%From $4_1$ to $6_4$ their htype is $[0,2]$ (two entangled handlebodies of genus $2$),
%htype is $[1,0,1]$ from $6_5$ to $6_7$, then $6_8$ to $6_{10}$ have htype $[1,2]$ and finally $6_{11}$ has
%htype $[2,0,1]$. The only entry that is not a $2$-sum is $6_1$.
}
\label{tab:hlinks}
\tabentryla{4}{a}{1}%
\tabentryla{6}{b}{1}%
\tabentryla{6}{c}{2}%
\tabentryla{6}{d}{3}%
\\
\smallskip
\tabentryla{6}{e}{4}%
\tabentrylb{6}{a}{5}%
\tabentrylb{6}{b}{6}%
\tabentrylb{6}{c}{7}%
\\
\smallskip
\tabentrylc{6}{a}{8}%
\tabentrylc{6}{b}{9}%
\tabentrylc{6}{c}{10}%
\tabentryld{6}{a}{11}%
\end{table}
%
%

%
%%%%%%%%%%%%%%%%%%%%%%%%%%%%%%%%%%%%%%%%%%%%%%%%%%%%%%%%%%%%%%%%%%%%%%%%%%%%%%
\section{Minimal diagrams; $1$-/$2$-sum}\label{sec:diagram}
 
%\subsection{Diagrams}  
\nada{
Given a (possibly nonconnected) spatial graph $\Gamma$ and $p\in \sphere-\Gamma$, 
we say a projection $\pi:\sphere-\{p\}\simeq \mathbb{R}^3\rightarrow \mathbb{R}^2$ is \emph{regular} with respect to $\Gamma$ if $\pi$ is in general position \cite{Zee:63} (see \cite[Section $4.4$]{Buo:03}). In particular, the image $\pi(\Gamma)$
is a plane graph with vertices the union of 
double points of $\pi$ and vertices of $\Gamma$. 

Given a regular projection $\pi$ with respect to $\Gamma$, a \emph{diagram} of $\Gamma$ is given by replacing each double point in $\pi(\Gamma)$
with a \emph{crossing}, where the broken arc corresponds
to the edge of $\Gamma$ underneath the edge the unbroken arc corresponds to in $\mathbb{R}^3\simeq \sphere-\{p\}$; see Fig.\ \ref{fig:over_under_crossings}.
The plane graph $\pi(\Gamma)$ is called \emph{the underlying graph} of the diagram.

A \emph{diagram} of a handlebody \knotlink\ $V$ is a diagram of one of its spines.
The \emph{crossing number} of a diagram of a spatial graph $\Gamma$ 
(resp.\ handlebody \knotlink\ $V$) is the number of the crossings the diagram has. The \emph{crossing number} $c(\Gamma)$ of $\Gamma$ 
(resp.\ $c(V)$ of $V$) 
is the minimum of the crossing numbers over all diagrams of $\Gamma$ (resp.\ $V$).
} 

\subsection{Handlebody-links}\label{subsec:diagram_handlebody_links}
We recall here terminologies pertaining minimal diagrams as our approach is based on detailed analysis on these diagrams. 
A minimal diagram of a handlebody-link $V$ is \emph{trivalent} if it represents a trivalent spine of $V$. Every handlebody-link with $\chi\leq -1$ 
admits a trivalent minimal diagram; see Fig.\ \ref{fig:turn_trivalent}. Our enumeration is done by first constructing a set $\mathcal{D}$ of trivalent minimal diagrams so that every irreducible handlebody-link with $\chi=-2$ is represented by some diagrams in $\mathcal{D}$, 
and then group diagrams representing equivalent handlebody-links. 
Hereinafter, minimal diagrams are assumed to be trivalent.

By \cite{Ish:08}, two minimal diagrams $D_1,D_2$ represent equivalent handlebody-links if they are related by a finite sequence of
generalized Reidemeister moves in Fig. \ref{fig:reidemeister_a} and \ref{fig:reidemeister_b} or some IH-moves in Fig. \ref{fig:IH-move}; in this case, we say $D_1,D_2$ are \emph{equivalent}. 

\subsection{Spatial graphs}\label{subsec:diagram_spatial_graphs}
\begin{figure}[ht]
\includegraphics[height=0.18\textwidth]{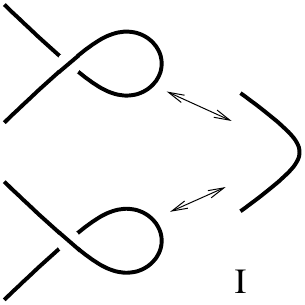}~~
\includegraphics[height=0.18\textwidth]{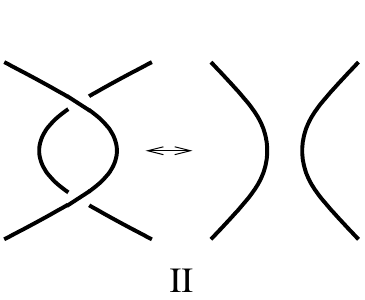}~~
\includegraphics[height=0.18\textwidth]{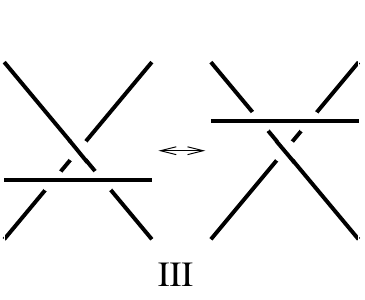}
\caption{Classical Reidemeister moves of type I, II, III.}
\label{fig:reidemeister_a}
\end{figure}

\begin{figure}[ht]
\includegraphics[height=0.18\textwidth]{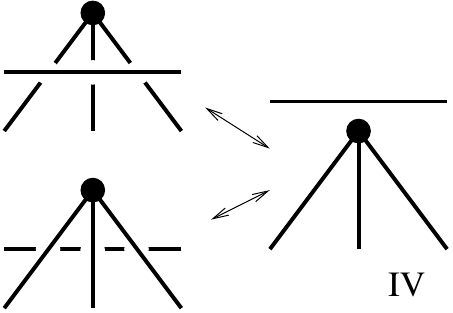}~~~~
\includegraphics[height=0.18\textwidth]{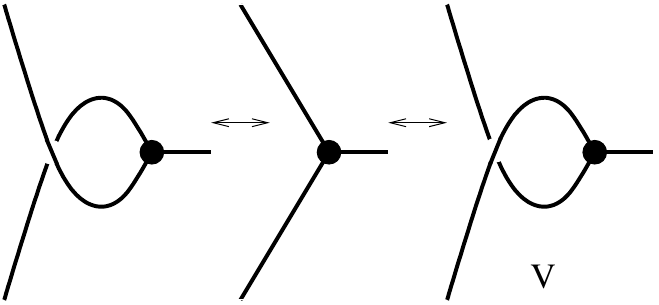}
\caption{Reidemeister moves IV and V involving a trivalent vertex.}
\label{fig:reidemeister_b}
\end{figure}

Our enumeration makes use of certain constructions on spatial graphs. Recall that, by \cite{Kau:89}, \cite{Yam:89}, two diagrams $D_1,D_2$ represent equivalent spatial graphs if 
they are related by a finite sequence of
generalized Reidemeister moves in Fig. \ref{fig:reidemeister_a} and \ref{fig:reidemeister_b}. In this case, we say $D_1,D_2$ are \emph{$R$-equivalent}. A diagram $D$ is \emph{$R$-minimal} if no finite sequence of
generalized Reidemeister moves in Fig. \ref{fig:reidemeister_a} and \ref{fig:reidemeister_b} can reduce its crossings.
In particular, $D$ is a minimal diagram 
of the spatial graph it represents.

A diagram $D$ of a spatial graph or handlebody-link has connectivity $k$ if $k$ is the minimum number of arcs whose removal disconnects the underlying graph
of $D$. The existence of trivalent vertices entails that $k \leq 3$. A diagram is \emph{split} 
if it has connectivity $0$.

\subsection{$2$-sum}\label{subsec:diagram_two_sum}
An \emph{order-$2$ vertex connected sum} of two spatial graphs $\Gamma_1,\Gamma_2$
is constructed as follows \cite{Mor:09a}:
First choose an edge $e_i$ in $\Gamma_i$, $i=1,2$. Secondly, consider a $3$-ball $B_i$ with $B_i\cap\Gamma_i=B_i\cap e_i$. Next, orient $e_i$ so that $e_i$ is going out of $B_i$ at a point $a_i\in e_i\cap \partial B_i$ and going in $B_i$ at a point $b_i\in e_i\cap \partial B_i$.
Then the \emph{order-$2$ vertex connected sum}  %abbreviated to \emph{$2$-sum}, 
is the spatial graph given by gluing the exteriors
$\Compl{B_1}$ and $\Compl{B_2}$ with an orientation-reversing homeomorphism
$h:\partial \Compl{B_1}\rightarrow \partial \Compl{B_2}$ with $h(a_1)=b_2$, $h(b_1)=a_2$, and the spatial graphs $\Gamma_1$, $\Gamma_2$ are called \emph{summands}. An order-$2$ vertex connected sum of $\Gamma_1,\Gamma_2$ depends on the choice of oriented arcs $e_1,e_2$ and is often not unique.
A handlebody-link equivalent to a regular neighborhood of an order-2 connected sum is called a \emph{$2$-sum}. We use $\Gamma_1\twosum \Gamma_2$ to denote the set of all $2$-sums obtained by $\Gamma_1,\Gamma_2$.

Let $V_1,V_2$ be two handlebody-links. Then a \emph{$1$-sum} $V_1\onesum V_2$ is constructed as follows:
First, choose a component $W_i\subset V_i$
and a $3$-ball $B_i$, $i=1,2$, such that $D_i:=V_i\cap B_i=\partial V_i\cap \partial B_i$ is a disk.
Then consider an orientation-reversing homeomorphism $h:\partial \Compl{B_1}\rightarrow \partial \Compl{B_2}$ with
$h(D_1)=D_2$. 
The resulting handlebody-link $(\sphere,V_1)\cup_h (\sphere, V_2)$ is called a \emph{$1$-sum} of $V_1$ and $V_2$. When $V_1$ or $V_2$ has more than one component, their $1$-sum may not be unique. 
We denote by $V_1\onesum V_2$ the set of all $1$-sums of $V_1,V_2$. 

\begin{remark}
A diagram with connectivity $k$, $k=1,2$, clearly represents a $k$-sum, but it is an open problem as to whether every $k$-sum admits a minimal diagram that has connectivity $k$. 

\end{remark}

\begin{figure}[ht]
\includegraphics[height=0.18\textwidth]{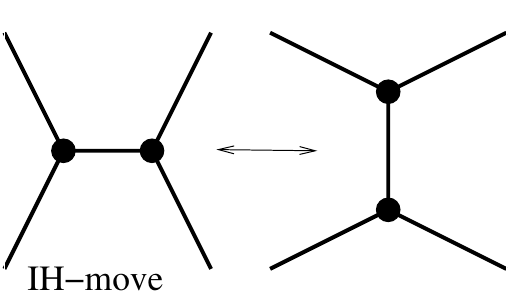}
\caption{IH-move.}
\label{fig:IH-move}
\end{figure}
\section{Completeness}\label{sec:completeness}
\begin{figure}[ht]
\def\svgwidth{0.9\columnwidth}
%% Creator: Inkscape 1.0beta2 (unknown), www.inkscape.org
%% PDF/EPS/PS + LaTeX output extension by Johan Engelen, 2010
%% Accompanies image file 'decomposing_type_2_diagram.pdf' (pdf, eps, ps)
%%
%% To include the image in your LaTeX document, write
%%   \input{<filename>.pdf_tex}
%%  instead of
%%   \includegraphics{<filename>.pdf}
%% To scale the image, write
%%   \def\svgwidth{<desired width>}
%%   \input{<filename>.pdf_tex}
%%  instead of
%%   \includegraphics[width=<desired width>]{<filename>.pdf}
%%
%% Images with a different path to the parent latex file can
%% be accessed with the `import' package (which may need to be
%% installed) using
%%   \usepackage{import}
%% in the preamble, and then including the image with
%%   \import{<path to file>}{<filename>.pdf_tex}
%% Alternatively, one can specify
%%   \graphicspath{{<path to file>/}}
%% 
%% For more information, please see info/svg-inkscape on CTAN:
%%   http://tug.ctan.org/tex-archive/info/svg-inkscape
%%
\begingroup%
  \makeatletter%
  \providecommand\color[2][]{%
    \errmessage{(Inkscape) Color is used for the text in Inkscape, but the package 'color.sty' is not loaded}%
    \renewcommand\color[2][]{}%
  }%
  \providecommand\transparent[1]{%
    \errmessage{(Inkscape) Transparency is used (non-zero) for the text in Inkscape, but the package 'transparent.sty' is not loaded}%
    \renewcommand\transparent[1]{}%
  }%
  \providecommand\rotatebox[2]{#2}%
  \newcommand*\fsize{\dimexpr\f@size pt\relax}%
  \newcommand*\lineheight[1]{\fontsize{\fsize}{#1\fsize}\selectfont}%
  \ifx\svgwidth\undefined%
    \setlength{\unitlength}{3685.03937008bp}%
    \ifx\svgscale\undefined%
      \relax%
    \else%
      \setlength{\unitlength}{\unitlength * \real{\svgscale}}%
    \fi%
  \else%
    \setlength{\unitlength}{\svgwidth}%
  \fi%
  \global\let\svgwidth\undefined%
  \global\let\svgscale\undefined%
  \makeatother%
  \begin{picture}(1,0.24615385)%
    \lineheight{1}%
    \setlength\tabcolsep{0pt}%
    \put(0,0){\includegraphics[width=\unitlength,page=1]{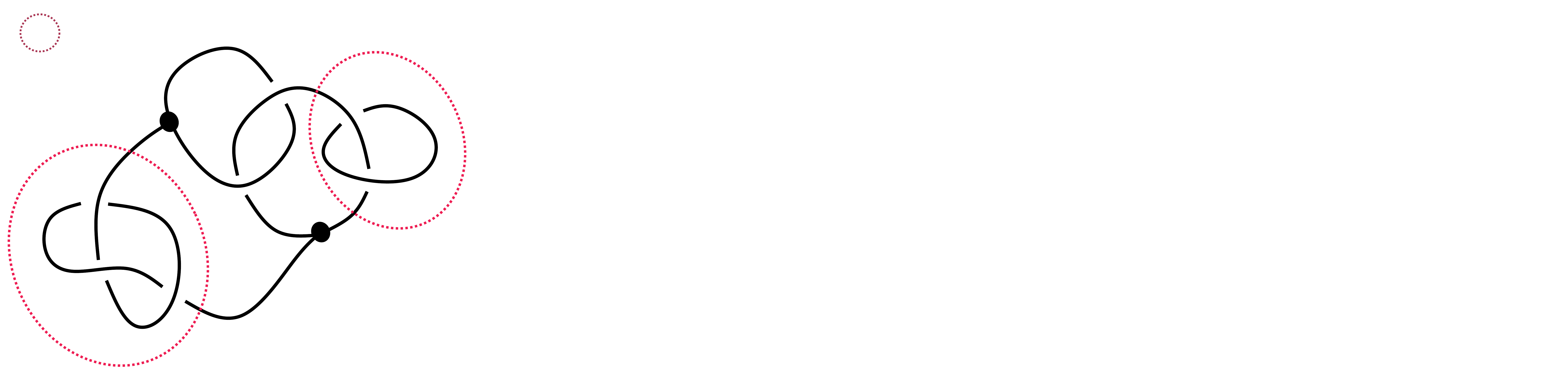}}%
    \put(0.0610321,0.22398191){\color[rgb]{0,0,0}\makebox(0,0)[lt]{\lineheight{1.25}\smash{\begin{tabular}[t]{l}{\footnotesize decomposing circle}\end{tabular}}}}%
    \put(0,0){\includegraphics[width=\unitlength,page=2]{decomposing_type_2_diagram.pdf}}%
    \put(0.70746202,0.01089694){\color[rgb]{0,0,0}\makebox(0,0)[lt]{\lineheight{1.25}\smash{\begin{tabular}[t]{l}{\footnotesize knot sum of spatial graphs}\end{tabular}}}}%
    \put(0,0){\includegraphics[width=\unitlength,page=3]{decomposing_type_2_diagram.pdf}}%
    \put(0.32583367,0.00734528){\color[rgb]{0,0,0}\makebox(0,0)[lt]{\lineheight{1.25}\smash{\begin{tabular}[t]{l}{\footnotesize induced spatial graphs}\end{tabular}}}}%
    \put(0,0){\includegraphics[width=\unitlength,page=4]{decomposing_type_2_diagram.pdf}}%
  \end{picture}%
\endgroup%

\caption{Decompose a $2$-connected, minimal diagram.}
\label{fig:decomposition_type_2}
\end{figure}

This section discusses completeness
of Tables \ref{tab:hknots} and \ref{tab:hlinks}.
Our aim is to build a set $\Dcal$ containing minimal diagrams of all non-split handlebody-links, both irreducible and reducible ones, up to six
crossings and with $\chi=-2$, up to equivalence and mirror image.
In other words, for each handlebody-link, $\Dcal$ contains at least one minimal diagram representing it.
%, but it can contain more than one  
%More than one diagram per handlebody is allowed, but at least one must be present.
The set $\Dcal$ will serve as the starting point for applying the \textit{dynamic process} described in
Section \ref{sec:uniqueness} and a correct construction of 
this set is essential to ensure completeness of 
Tables \ref{tab:hknots} and \ref{tab:hlinks}.

We shall treat differently diagrams with connectivity $k$ with $k = 1,2,3$%
\footnote{$k=0$ is excluded since it would lead to split handlebody-links.}
in order to obtain
$\Dcal$ as 
\begin{equation}\label{eq:D}
\Dcal = \Done \cup \Dtwo \cup \Dthree;
\end{equation}
the case $k=1$, producing $\Done$ and treated
in Subsection \ref{sec:typeone}, yields handlebody-links that are evidently reducible.
The converse may not hold: it is possible that a reducible handlebody-link possesses
no \textit{minimal} diagram with connectivity $1$\footnote{%
A diagram with connectivity $1$ must exist, as clear from the definition of reducibility, 
%(Definition \ref{def:R-minimal}),
however it is possible that such diagrams all have crossing numbers 
larger than the crossing number of
the represented handlebody-link.
No such examples have been found, and it is, in fact, conjectured to 
never occur. This is, nevertheless,
still open.%
}.
For this reason we devote Section \ref{sec:irreducibility} to prove that all entries in Tables \ref{tab:hknots}
and \ref{tab:hlinks} are indeed irreducible.

Having ruled out handlebody-links possessing a minimal diagram with connectivity $1$
we can consider diagrams with connectivity $2$,
representing handlebody-links that are the result of a $2$-sum (see Subsection \ref{subsec:diagram_two_sum}).
This is done in Subsections \ref{sec:typetwo} and \ref{sec:typetwospecial} where we build all possible results of
$2$-sums involving tangles with connectivity $3$ (representing a spatial graph) as factors and we can avoid considering resulting
diagrams that are shown to be equivalent to a diagram with connectivity $1$.
It is possible to have diagrams in $\Dtwo$ that are equivalent to some diagram with connectivity $1$; this is however harmless.

Finally, 
in Subsection \ref{sec:typethree}
we construct $\Dthree$ 
 by considering only diagrams with connectivity $3$.
Again it is possible, but harmless, for some of them to 
be equivalent to a diagram with connectivity $k$, $k \in \{1,2\}$.
Set $\Dthree$ is the most critical one to construct, since its entries are not built from simpler objects and yet
we need to be sure that \emph{all} possible handlebodies (up to equivalence and mirror image and up to six crossings)
are represented, 
at least those that cannot be obtained as $2$-sums or $1$-sums.
They are obtained using a software code developed by the 
second author that constructs all possible plane graphs with
four trivalent vertices and $n \in \{0,\dots,6\}$ quadrivalent 
vertices.

In all three cases ($k = 1,2,3$) we are allowed to drop all diagrams that are found to be nonminimal or, in case
of a set of equivalent diagrams, to retain only one of them.

We shall occasionally refer to files like \texttt{embeddings\_*},
built using the above-mentioned computer software.
They can be found in \cite{bigtable:web}.

%Recall first that a minimal diagram of a 
%non-split, irreducible handlebody-\knotlink\ has
%either connectivity $2$ or $3$ (see Definition \ref{def:diagram_connectedness}).
%minimal diagrams with connectivity $3$ 
%are obtained from a software code,
%and handlebody \knotlinks\ represented by 
%minimal diagrams with connectivity $2$
%are recovered by knot sum of spatial graphs. 
 
\subsection{Minimal diagrams with connectivity 3}
\label{sec:typethree}

\begin{figure}[h]
\includegraphics[height=0.1\textwidth]{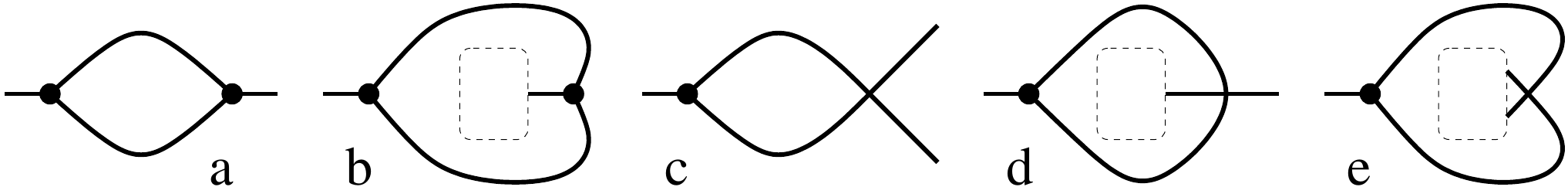}
\includegraphics[height=0.1\textwidth]{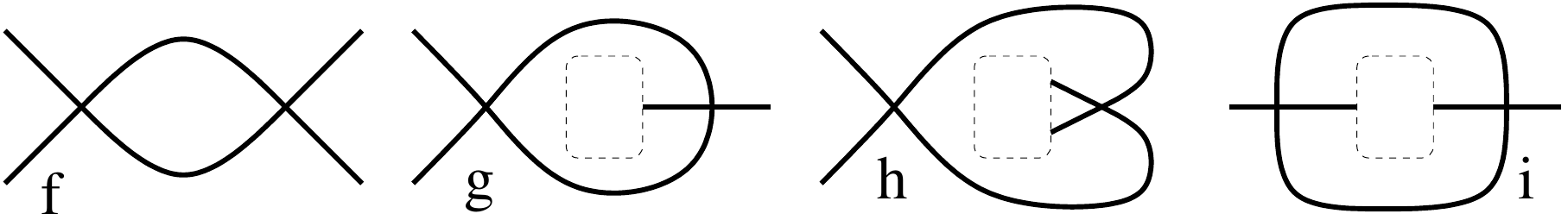}
\caption{Possible configurations for loops and double arcs.}
\label{fig:loops_double_arcs}
\end{figure}

Our goal here is to find all handle\-body-links with
up to six crossings that have a $3$-connected minimal diagram.
To do so, we first enumerate all $3$-arc-connected plane graphs with 
four trivalent vertices and
up to six quadrivalent vertices,
and then replace each quadrivalent vertex with a crossing.%
\footnote{
Note that the Whitney theorem \cite{Whi:32} is not 
applicable here since a $3$-connected graph may not be $3$-vertex-connected.%
\draftMP{[example: embeddings \#27 and \#28 with $4$ quadrivalent vertices (and $4$ trivalent vertices)}}
\begin{center}
\begin{tabular}{c|rrrrrrr}
$n$ & 0 & 1 & 2 & 3 & 4 & 5 & 6 \\
\hline
Number of plane graphs & 1 & 1 & 4 & 18 & 110 & 677 & 4752
\end{tabular}
\end{center}
The first part is done by a code devised by the second author, resulting in
$5563$ such plane graphs, $4752$ of which with six quadrivalent vertices,
up to mirror image.
This enumeration is obtained by first obtaining all abstract graphs with the given
number of trivalent and quadrivalent vertices (also allowing for double arcs between
two nodes, since bigons are naturally present in a typical diagram), then checking whether
each graph is planar and, if yes, all distinct embeddings that match the following
conditions are recorded:
\begin{enumerate}
\item
each of them has edge-connectivity $3$ as an abstract graph, i.e. all resulting diagrams
must be $3$-connected;
\item
their double arcs can only connect two quadrivalent vertices as abstract graphs;
\item
their double arcs only form a ``bigon'' in the planar embedding (a polygon with two sides corresponding
to case `f' in Fig.\ \ref{fig:loops_double_arcs});
\item two vertices cannot be joined by more than two arcs.
\end{enumerate}
Condition (3) in not actually 
enforced by the code, since it follows from the other requirements.
More specifically cases ``a, b, c, d, e'' are ruled out by condition (2) whereas cases
``g, h, i'' lead to diagrams with connectivity $1$ or $2$ contradicting condition (1).\footnote{%
Condition $(1)$ is compatible with ``i'' if we allow the four dangling arcs to connect
to each other,
however this is against the requirement that we have trivalent vertices.
In case of no trivalent nodes (i.e. the case of standard knots/links) we contradict condition (4).}

The reason why we can filter the plane graphs by 
enforcing conditions (2),(3),(4) is because otherwise
we either get a non-R-minimal diagram, or we contradict $3$-connectedness
(see Fig.\ \ref{fig:loops_double_arcs},
where ``a, b, c, d, e'' illustrate those double arcs 
connecting at least one trivalent vertex and ``g, h, i''  
those connecting two quadrivalent 
vertices with a non-bigon configuration.)
Note that this is true by recalling that we have (at least three) trivalent vertex, otherwise we should
also specifically allow for the trivial theta-graph and for the Hopf link $\op{L2a1}$.

\begin{figure}
\begin{subfigure}{.48\linewidth}
\centering
\includegraphics[scale=.12]{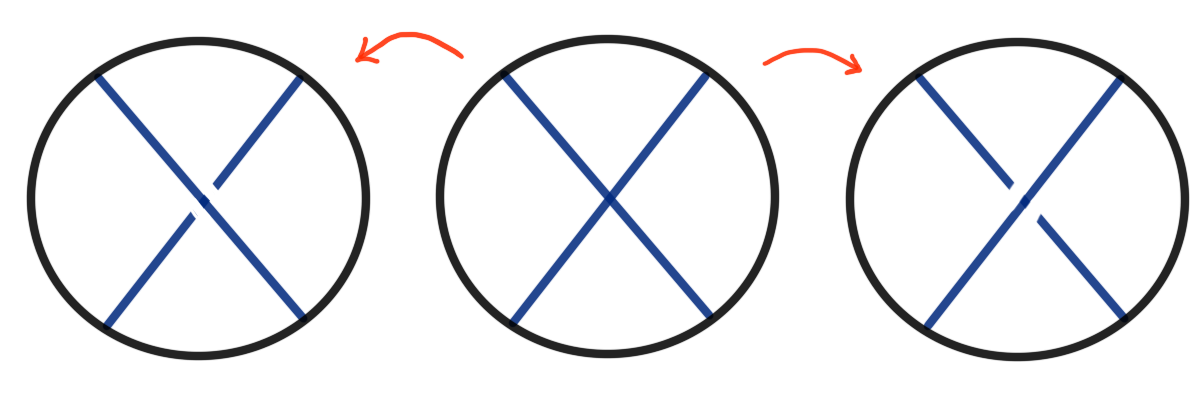}
\caption{Turn $4$-valent vertices into crossings.}
\label{fig:over_under_crossings}
\end{subfigure}
\begin{subfigure}{.48\linewidth}
\centering
\includegraphics[scale=.12]{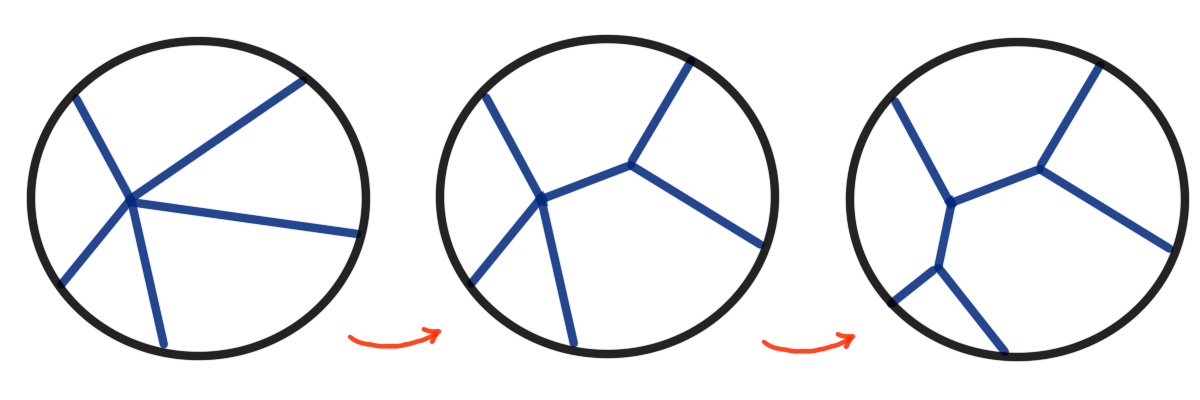}
\caption{Turn vertices into trivalent ones.}
\label{fig:turn_trivalent}
\end{subfigure}
\caption{}
\end{figure}

\begin{figure}
        \includegraphics[height=0.3\textwidth]{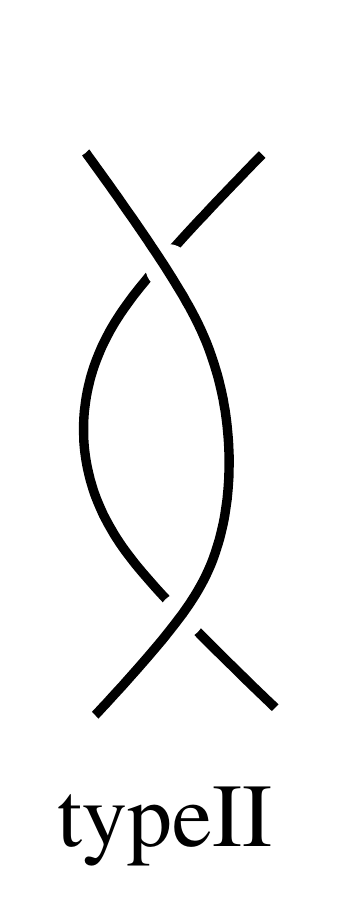}
        ~~
        \includegraphics[height=0.3\textwidth]{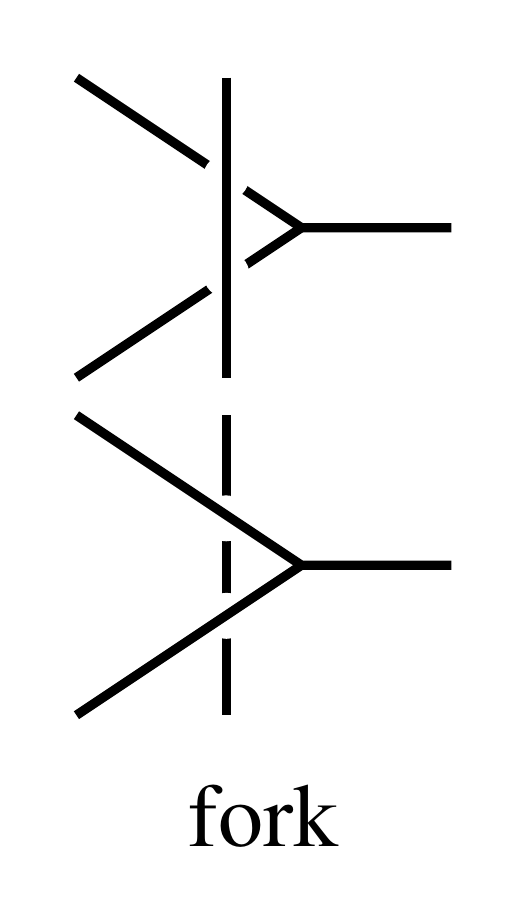}
        ~~
        \includegraphics[height=0.3\textwidth]{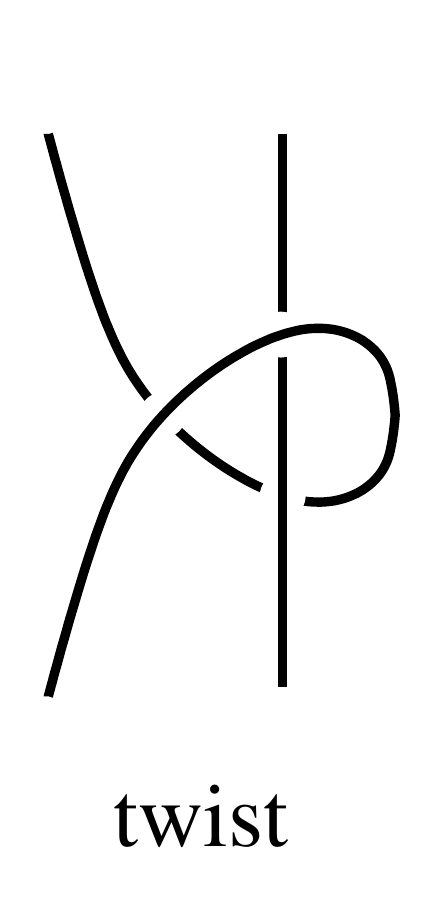}
        ~~
        \includegraphics[height=0.3\textwidth]{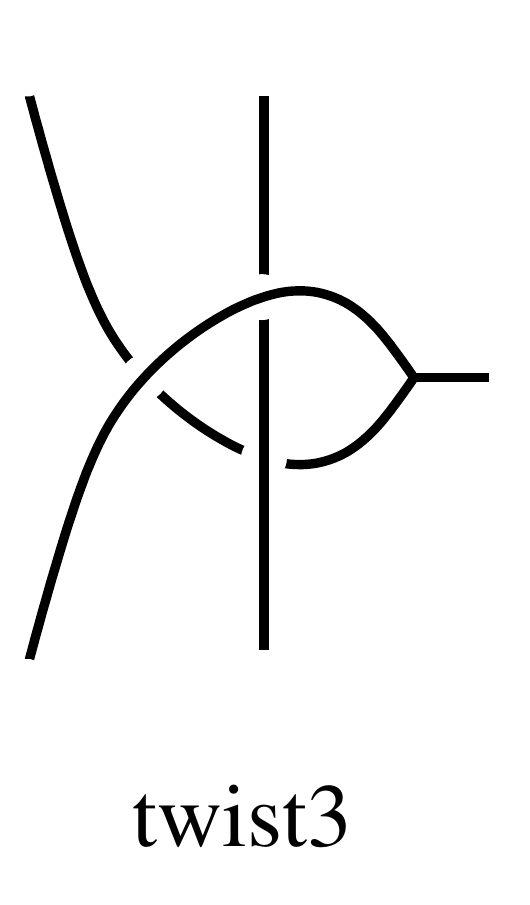}
        ~~
        \includegraphics[height=0.3\textwidth]{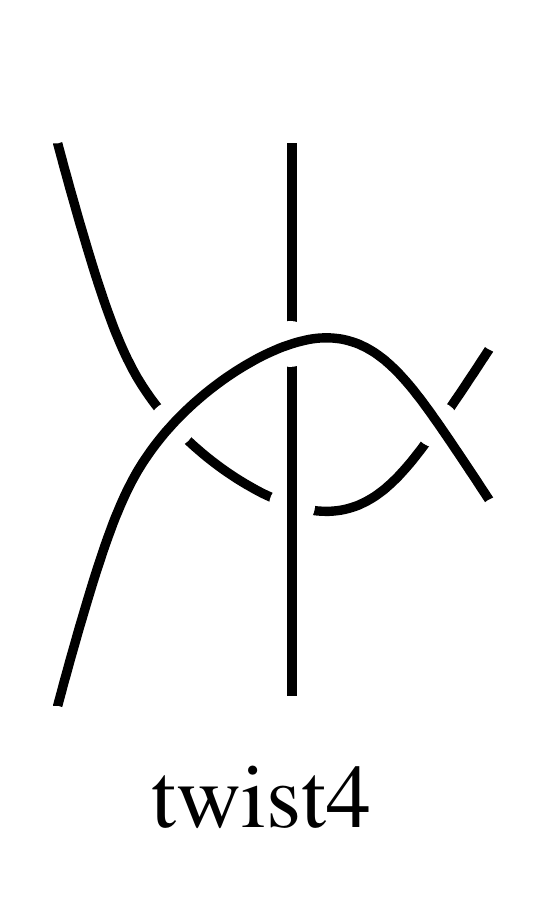}
        \caption{Crossing-reducing moves.  They consist of a sequence of Reidemeister moves
        as follows (from left to right):
        II~; II$^{-1}$ IV~; III I~; III V~; III II~. Notation ``II$^{-1}$" means that move II
        is performed in the ``crossing-increasing'' direction.}
        \label{fig:autononminimal}
\end{figure}

\begin{figure}
        \includegraphics[height=0.3\textwidth]{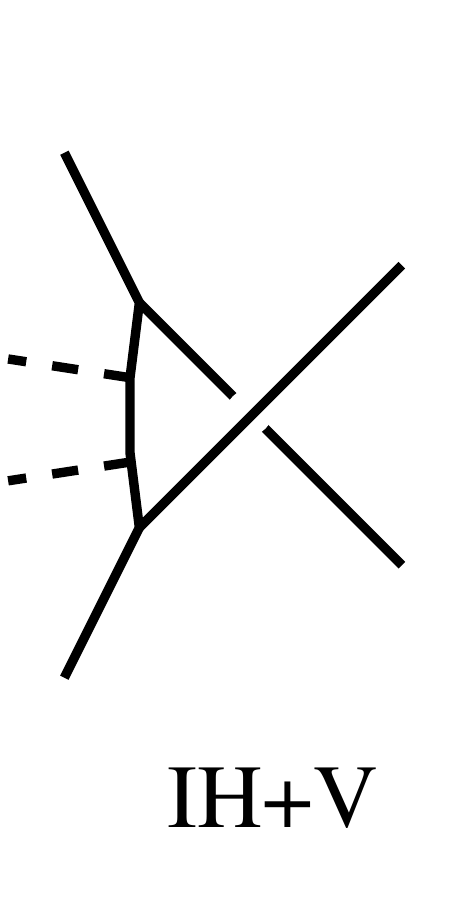}
        ~~~~~~~
        \includegraphics[height=0.3\textwidth]{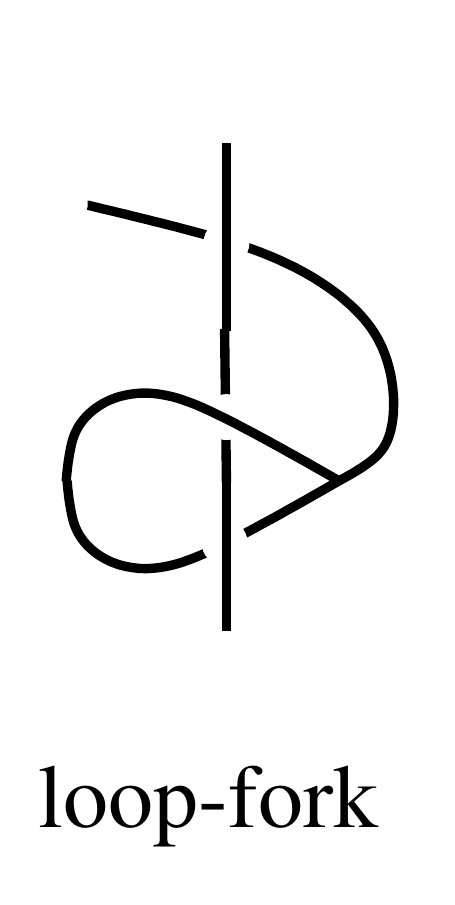}
        ~~~~~~~
        \includegraphics[height=0.3\textwidth]{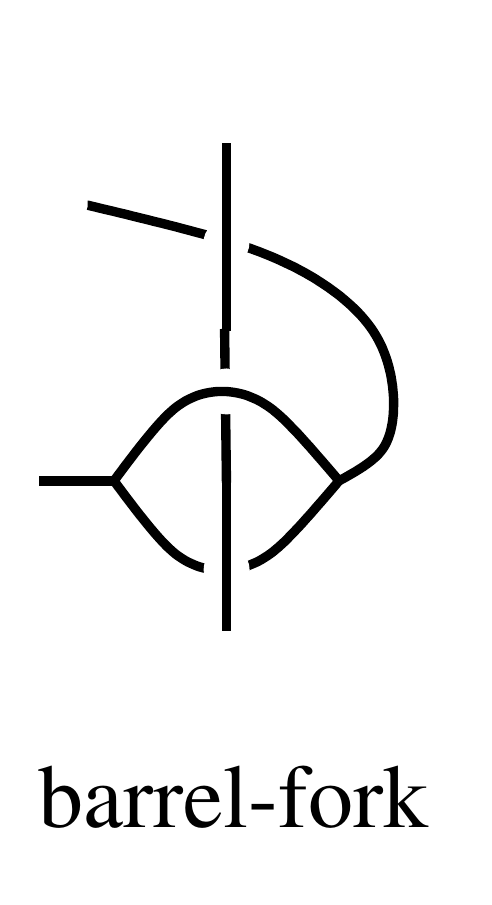}
        \caption{More crossing-reducing moves. ``IH+V'' is a sequence of $k-1$ IH-moves, where $k$ is the
	number of trivalent vertices involved, followed by a V
        move (Fig. \ref{fig:reidemeister_b});
        ``loop-fork'' and ``barrel-fork'' are performed by twisting the loop/barrel (see Fig.
        \ref{fig:autopairs}) and then applying the fork move
        (Fig.\ \ref{fig:autononminimal}).}
        \label{fig:moreautononminimal}
\end{figure}

\begin{figure}
        \includegraphics[height=0.3\textwidth]{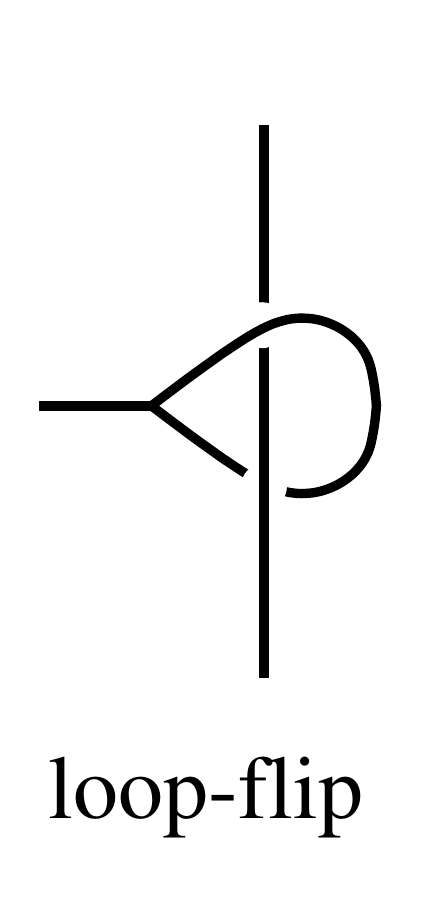}
        ~~~~~~ 
        \includegraphics[height=0.3\textwidth]{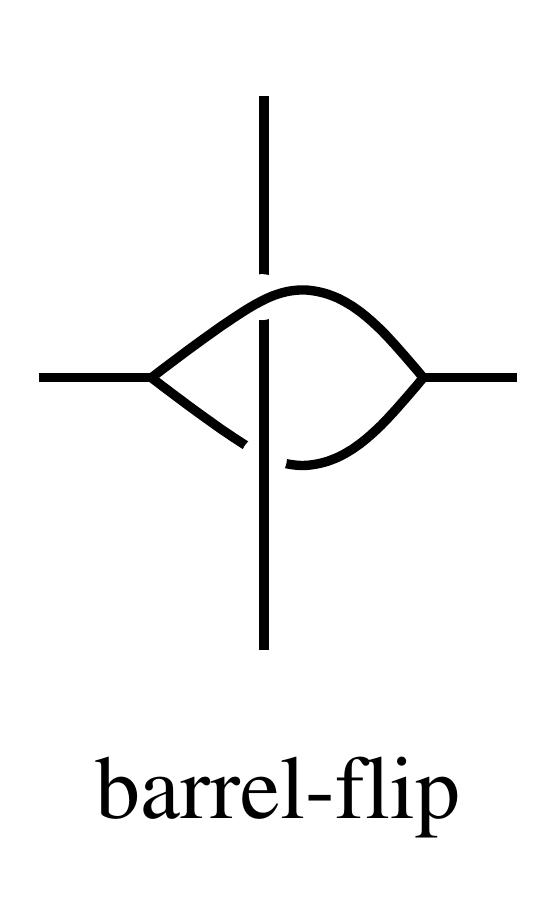}
        \caption{Loop/barrel flipping moves.  They consist of the following sequence of Reidemeister moves
        ``loop-flip'': II$^{-1}$ V III I resulting in a reversal of the two overcrossings.
        ``barrel-flip'': same sequence as for ``loop-flip'' with a V move replacing I,
        also resulting in a reversal of the overcrossings.}
        \label{fig:autopairs}
\end{figure}

Since, for each quadrivalent vertex, there are two possible choices of crossings (see Fig.\ \ref{fig:over_under_crossings}), these $5559$ plane graphs produce an unmanageable number,
$\approx 5559\cdot 2^5$, of diagrams, up to mirror image. To reduce the number of diagrams, a software developed by the third author is implemented to drop diagrams that
\begin{itemize}
\item{} are recognized to be nonminimal via moves in Figs.\  \ref{fig:autononminimal} and \ref{fig:moreautononminimal}.
%, i.e. easily transformed by using the generalized Reidemeister
%moves to a diagram with fewer crossings (see Subsection \ref{sec:autononminimal}),
\end{itemize}
In addition, one choice is selected for
any pair of over/under crossings
\begin{itemize}
        \item{} that result in embeddings having a portion depicted in Fig.\ \ref{fig:autopairs} or
        \item{} that are directly connected by an IH-move (Fig.\ \ref{fig:IH-move}).
\end{itemize}

The filtering process described in these three items does not guarantee against nonminimal
diagrams.
Those must be manually uncovered during the \textit{dynamic process} described in Section \ref{sec:dynamicapproach}.

Consequently, we obtain a refined set $\widehat{\Dthree}$ containing %
% $910$ (before reducing up to IH-move)
$8090$ diagrams ($7268$ with six crossings) where, for every $n$-crossing, 
$n \leq 6$,
handlebody-link $V$ with a $3$-connected minimal diagram, there is a minimal diagram in $\widehat{\Dthree}$ representing $V$.
We then remove entries from $\widehat{\Dthree}$ to
arrive at the final number $388$ of $\Dthree$ in two steps, first all situations depicted in Figures \ref{fig:autononminimal}
and \ref{fig:moreautononminimal} with the exception of the IH+V move are recognized by the \texttt{appcontour} software
and used to generate a set of 
$$
8090 = 0 + 0 + 2 + 11 + 96 + 713 + 7268$$
diagrams ($7268$ with six crossings).

This set is then reduced to 
$$
388 = 0 + 0 + 0 + 0 + 4 + 21 + 363
$$
entries, set $\Dthree$ ($363$ with six crossings) 
 by automatically recognizing situations matching the move ``IH+V'' of Figure
\ref{fig:moreautononminimal} 
and, for any set of diagrams that are directly connected by IH-moves, only one is retained.

% end \nada

% end \nada

\begin{figure}
\includegraphics[height=0.20\textwidth]{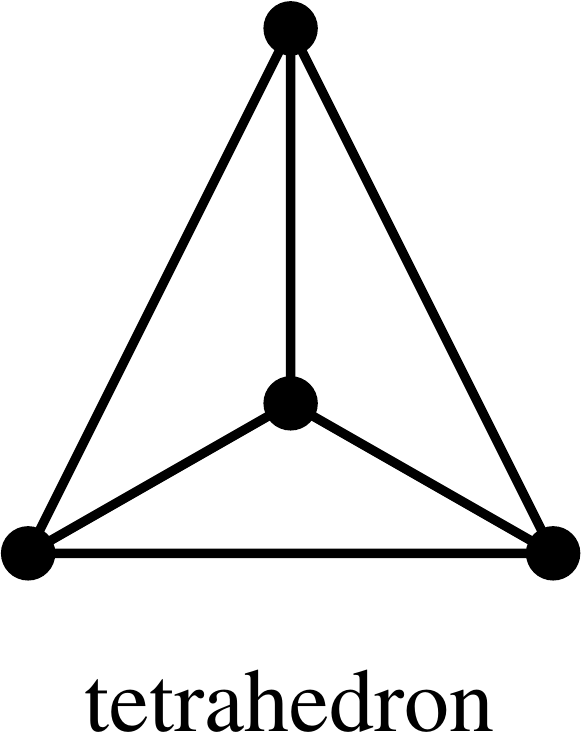}~~
\includegraphics[height=0.20\textwidth]{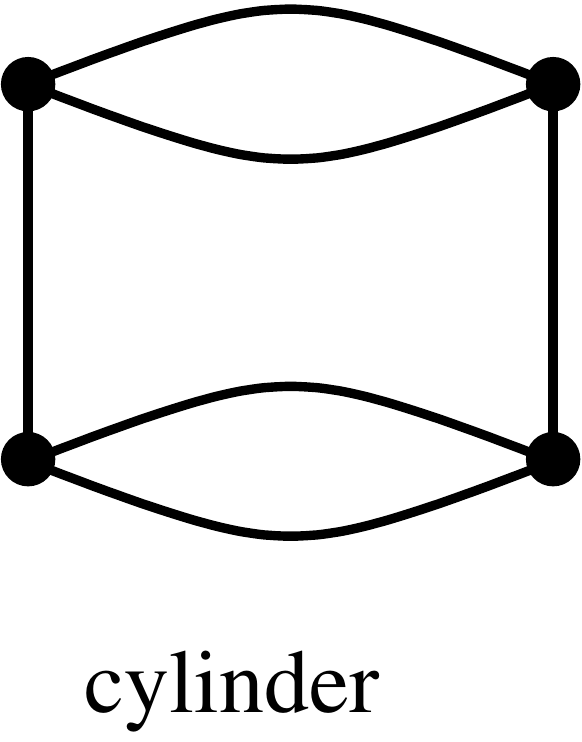}~~
\includegraphics[height=0.20\textwidth]{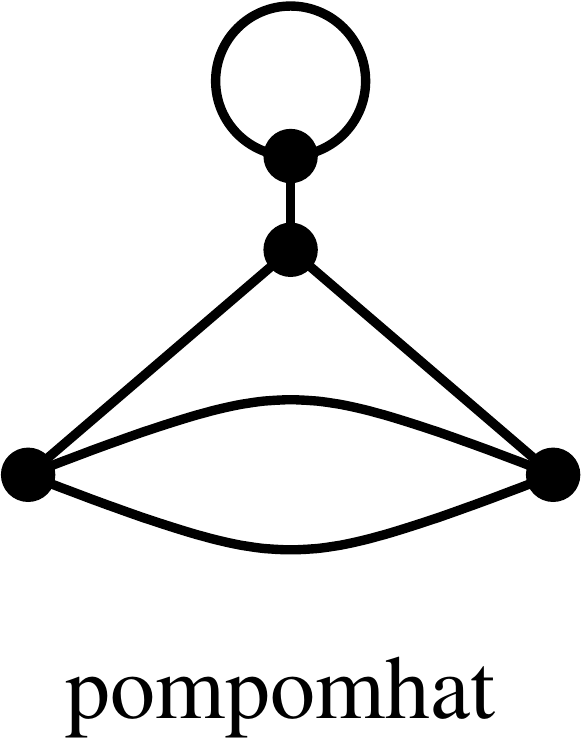}~~
\includegraphics[height=0.20\textwidth]{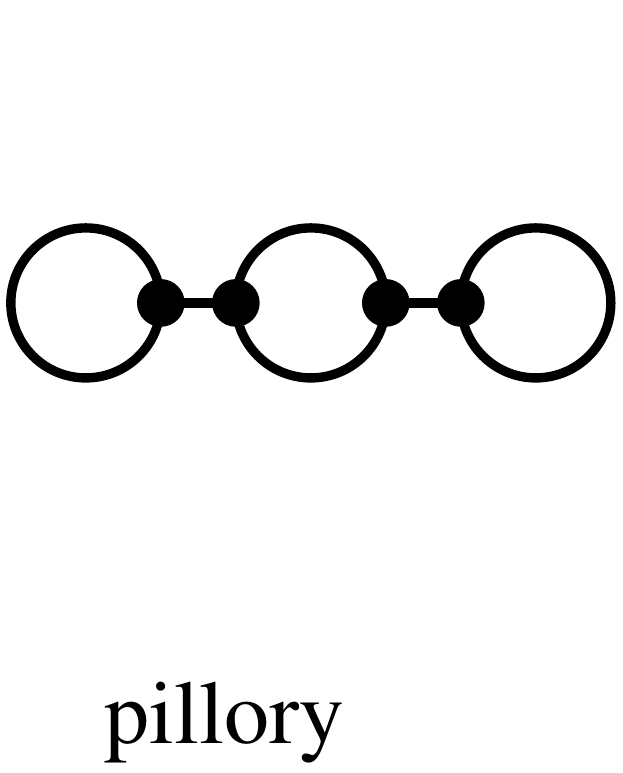}~~
\includegraphics[height=0.20\textwidth]{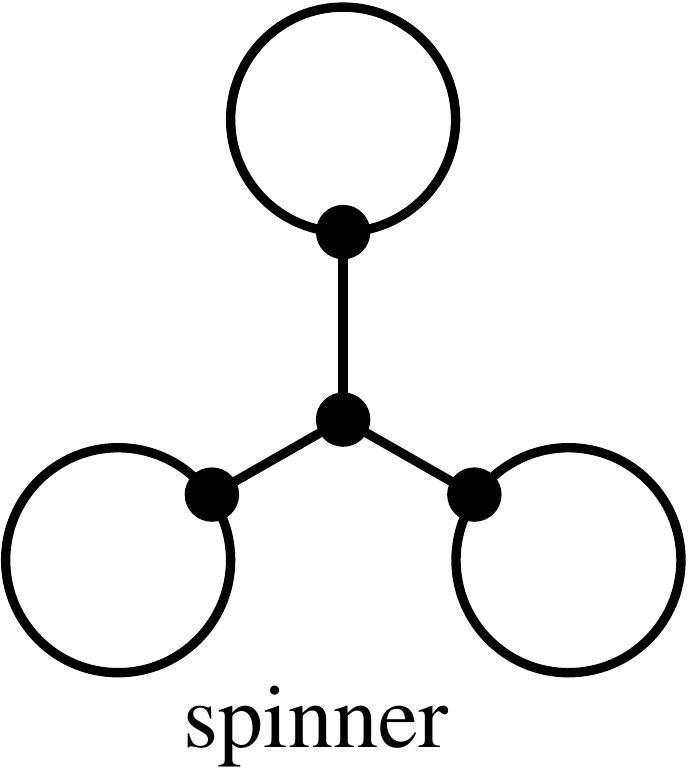}
\caption{The five trivial spatial graphs with fancy names (\textit{gtypes}).
They correspond
to the well-known five trivalent connected genus-3 graphs.
In the literature they are usually referred to as (000), (020), (111), (212), (303)
respectively,
see e.g. \cite[Figure 4]{BrJoMoSt:15}.
}
\label{fig:trivialgraphs}
\end{figure}

\subsection{Minimal diagrams with connectivity $2$}\label{sec:typetwo}
%%%%%copied from existence2
%%%%%%%%%%%%%%%%%%%%%%%%%%%%%%%%%%%%%%%%%%%%%%%%%%%%%%%%%%%%%%%%%%%%%%%%%%%%%%%
Recall that a diagram $D$ with $2$-conn\-ectivity 
can be decomposed into a \textit{maximal $2$-sum}, i.e. it can be
decomposed into
finitely many simpler (nontrivial) tangle diagrams such that each 
associated diagram of spatial graphs has connectivity $3$ or $4$  
(Fig.\ \ref{fig:decomposition_type_2}).
Furthermore, if $D$ is R-minimal (see Section 
\ref{subsec:diagram_spatial_graphs}) each induced spatial graph 
diagram is also R-minimal. In particular, a minimal diagram 
with connectivity $2$
can be recovered by performing the order-2 vertex connected sum
(in short $2$-sum) between spatial graphs admitting a R-minimal diagram 
with connectivity $k>2$.

\begin{figure}
\includegraphics[height=0.20\textwidth]{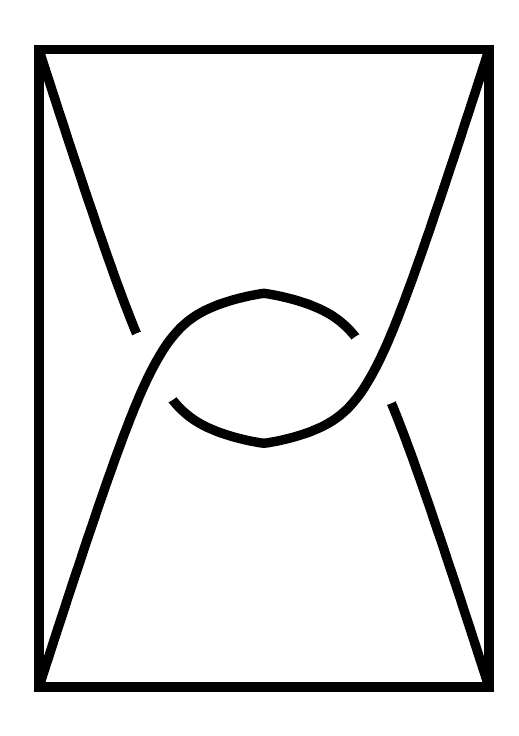}~~~~~~
\includegraphics[height=0.20\textwidth]{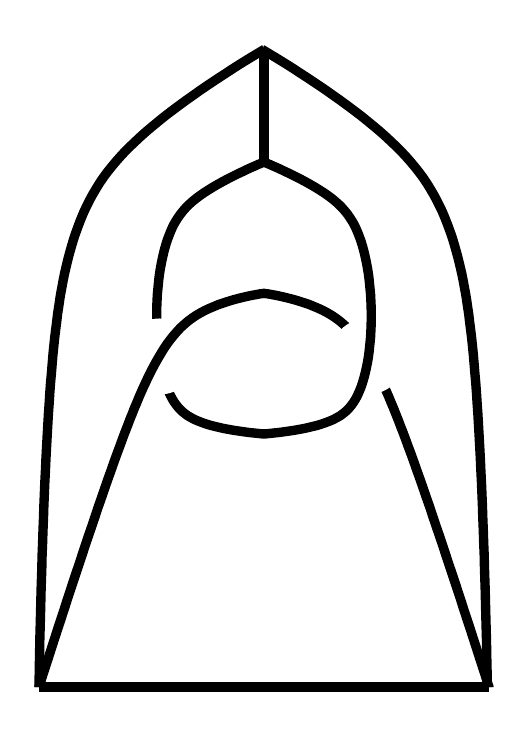}
\caption{R-minimal spatial graphs with arc-connectivity $3$ and $2$ crossings, as
uncovered by our analysis. Their gtype is \textit{cylinder} (left) and \textit{pompomhat} (right).}
\label{fig:x2-graphs}
\end{figure}

\begin{figure}
\includegraphics[height=0.20\textwidth]{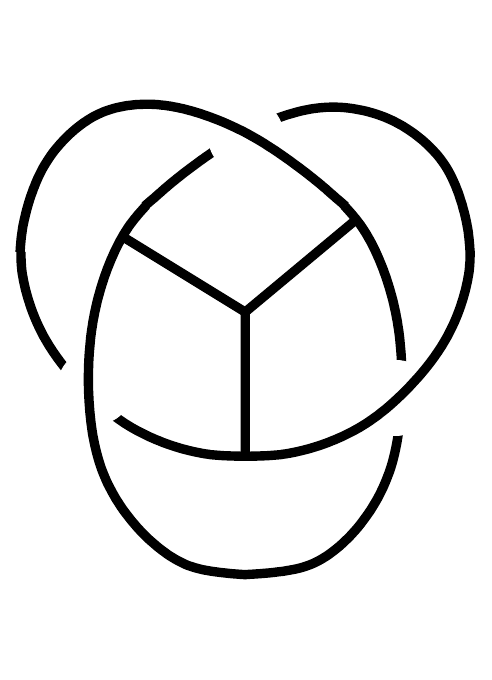}\!\!\#5~
\includegraphics[height=0.20\textwidth]{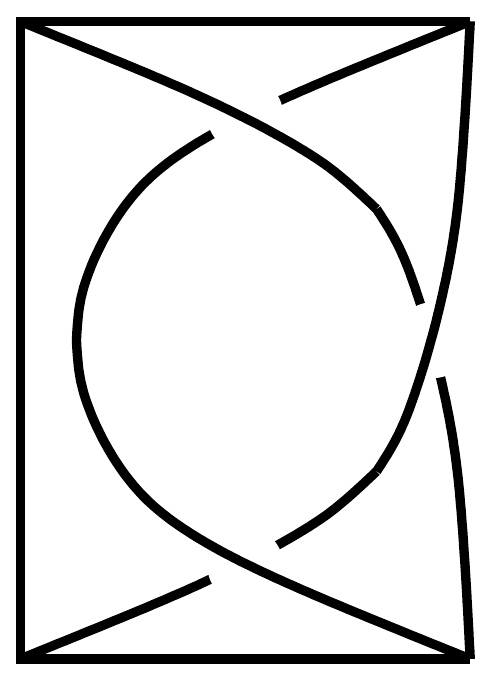}\#13~
\includegraphics[height=0.20\textwidth]{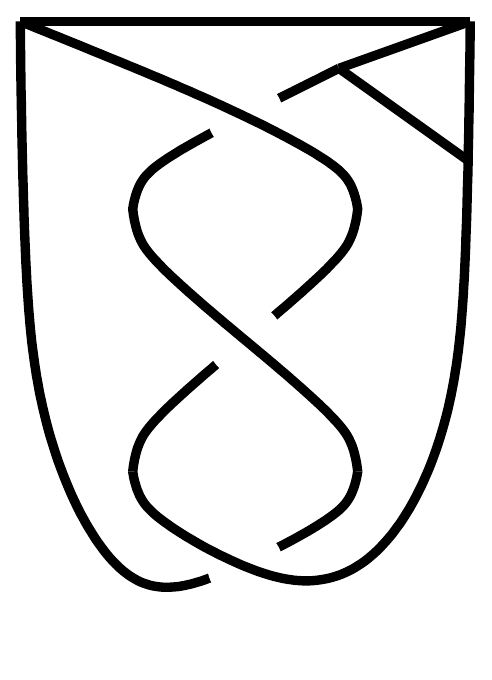}\!\!\!\!\#14~
\includegraphics[height=0.20\textwidth]{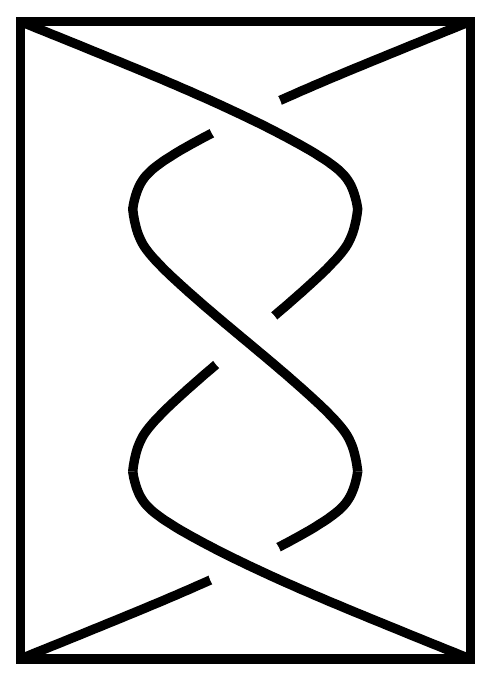}\#17~
\includegraphics[height=0.20\textwidth]{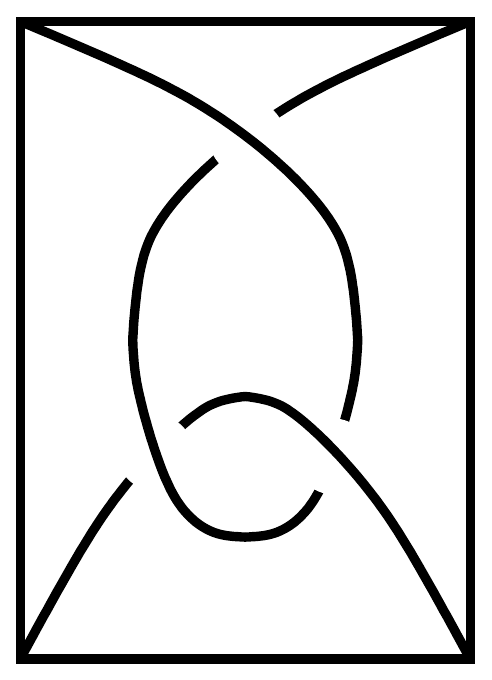}\#18
\caption{R-minimal spatial graphs with arc-connectivity $3$ and $3$ crossings
tagged with their reference in file \texttt{embeddings\_3}.}
\label{fig:x3-graphs}
\end{figure}

\begin{remark}
{\sl Chirality and reversibility.}
In general there are many distinct possible results of a $2$-sum, since it depends
on the choice of the arcs involved for each factor and a choice of their orientation.
To rule out the dependence on the orientation we observe that if either one of the
two arcs involved is {\it reversible} within its factor,
then the result does no longer depend on the choice
of orientation.
An arc of a spatial graph is {\it reversible} if there exists an orientation preserving
homeomorphism of $\sphere$ that maps the spatial graph onto itself while reverting
the orientation of the selected arc.
When performing a $2$-sum between two spatial graphs having $3$-connected minimal diagram with
a resulting diagram with up to six crossings then one of the two factors must have up to two
crossings with the only exception of a $2$-sum between two spatial graphs with $3$-connected
minimal diagram and exactly three crossings.
Now let us enumerate all spatial graphs with a $3$-connected minimal diagram and increasing number of
crossings.
\begin{itemize}
\item (0 crossings) We have the trivial theta-curve $\op{tG0}$ (genus $2$) and the trivial \texttt{tetrahedron}
$\op{T0}$ of Figure \ref{fig:trivialgraphs}.  They are both clearly reversible (meaning that all 
their arcs are reversible);
\item (1 crossing) There are no such spatial graphs;
\item (2 crossings) We have the Hopf link $\op{L2a1}$, $\op{hG2}$ (entry $2_1$ in \cite{Mor:07b}).
They are easily seen to be reversible.
Then we have the two spatial graphs shown in Figure \ref{fig:x2-graphs}, the first is reversible,
the second is presumably not reversible, however it can only be involved in a $2$-sum with a standard
knot/link with up to four crossings, all of which are reversible;
\item (3 crossings) We have the trefoil knot $\op{K3a1}$, $\op{tG3}$ (entry $3_1$ in \cite{Mor:09a})
(both reversible) and
the spatial graphs shown in Figure \ref{fig:x3-graphs}.  However these can only be involved in a $2$-sum
with a standard knot/link with up to three crossings, all of which are reversible.
\end{itemize}
In conclusion in all cases at least one of the two factors is reversible, so that we need not consider both
choices for the arcs orientation.
Concerning chirality, if we perform the $2$-sum of two chiral objects, the result can
be nonunique, even if considered up-to-mirror-image.
The only circumstance where this can happen in what follows is when both factors are
the trefoil knot $\op{K3}$ ($\op{3_1}$ in Rolfsen table) or the theta-curve $\op{tG3}$ 
(entry $3_1$ in \cite{Mor:09a})
Indeed we do have two distinct entries $6_{12}$ and $6_{13}$ in Table \ref{tab:hknots}.
This type of behaviour also occurs when enumerating reducible 
handlebody-links,
see Table \ref{tab:typeone} and Fig. \ref{fig:typeonedetail}.
\end{remark}

Our aim is to
find all possible maximal 2-sums among $n$ tangle diagrams with a
result having up to six crossings, $n \geq 2$.

{\bf Notation.} In the rest of this subsection we shall stick with the
following notation  for a $3$-connected handlebody link:
$$
{\rm genus}~ g+
$$
stands for a 
handlebody-link having genus
$g + \numcc - 1$, where $\numcc$ is the number of connected components of the handlebody; for example, genus $g$ for a handlebody-knot, and $g+1$ for a $2$-component handlebody-link, and so on.

Recalling \eqref{eq:h_type}, we shall denote by:
\begin{itemize}
\item $\fthree$
 the diagram of a 
generic $3$-connected
handlebody-link of genus $3+$, i.e. of type $[*,0,1]$ or type $[*,2]$,
\item
 $\ftwo$ the diagram of a generic handlebody-link of genus $2+$, i.e.
of type $[*,1]$.
\end{itemize}

For a generic $2$-sum of $n$ handlebody factors that results in a handlebody with
Euler characteristic $-2$ (i.e. genus $3+$) we can distinguish two cases.
\begin{itemize}
\item[a)]
One of the factors has genus $3+$ and the remaining factors are all
classical links %\knotslinks\ 
(genus $1+$);
\item[b)]
two of the factors have genus $2+$ and the remaining are all classical links.
%standard \knotslinks.
\end{itemize}

In any case at least $n-2$ of the factors must be the diagram of a classical link,
%standard \knotlink, 
each one
of them contributing for at least $2$ crossings (exactly
$2$ when the factor is the Hopf link $\op{L2a1}$)
for a total of $2(n-2)$ crossings.  Since we are interested in the case of a maximum of
$6$ crossings we can then 
conclude that 
$$
n \leq 5.
$$
We shall enumerate all possible constructions with up to $6$ crossings for values of $n$
from $5$ down to $n = 2$ and deal separately with the case $n = 2$, six crossings and one
factor being the Hopf link, i.e. $H = \fthree \twosum \op{L2a1}$ with $\fthree$ with $4$ crossings.

From now on, the symbol $\flink$ stands for a classical link.
%standard knots/links.

%\begin{remark}
%When dealing with handlebodies a set of $2$-sums can be described using an abstract tree where
%each node represents a factor (a genus $g+$ handlebody knot/link) and each arc represents a
%$2$-sum between the factors represented by the two nodes connected by the arc.
%If $n \geq 4$ this structure is no-longer faithfully representable as a linear expression using parentheses
%and the $2$-sum operator.
%Still, each $2$-sum can be applied in different ways and requires the indication of which arc of
%the two spatial graphs are involved and the choice of orientation for those arcs.
%\end{remark}

\begin{figure}
        \includegraphics[height=0.12\textwidth]{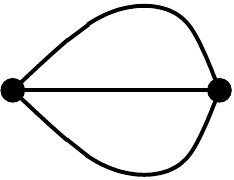}$\op{tG0}$
        \includegraphics[height=0.15\textwidth]{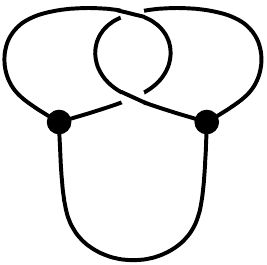}$\op{hG2}$
        \includegraphics[height=0.15\textwidth]{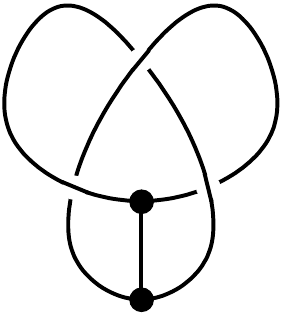}$\op{tG3}$
        \includegraphics[height=0.15\textwidth]{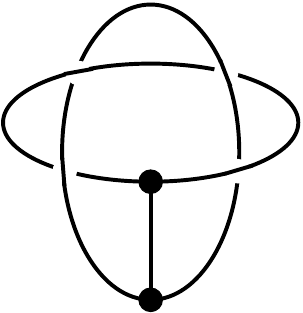}$\op{hG4}$\\
        \vspace*{1.5em}
        \includegraphics[height=0.15\textwidth]{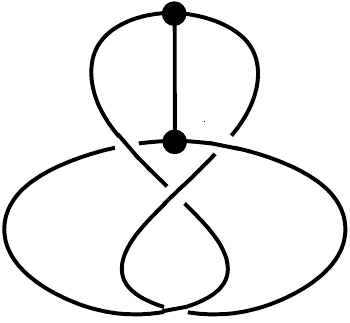}$\op{tG4}$
        \includegraphics[height=0.15\textwidth]{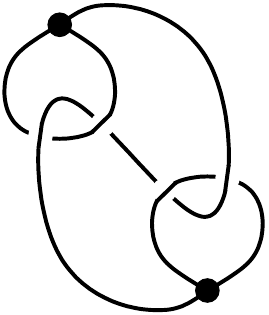}$\op{HK4}$
        \includegraphics[height=0.15\textwidth]{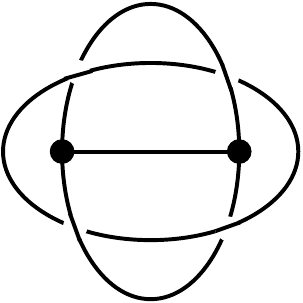}$\op{{\ell}G4_1}$
        \includegraphics[height=0.15\textwidth]{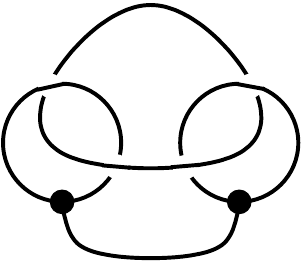}$\op{{\ell}G4_2}$
        \caption{Spatial graphs of genus $2+$ up to four crossings.
        $\op{hG2}$ and $\op{hG4}$ are handcuff graphs and correspond respectively to
        $2_1$ and $4_1$ of \cite[Table 1]{Mor:07b}; $\op{HK4}$ is a nonprime handcuff graph
        obtained as the $3$-sum of two copies of entry $2_1$ in \cite{Mor:07b};
        $\op{tG3}$ and $\op{tG4}$ are theta-curves corresponding respectively to
        $3_1$ and $4_1$ of \cite[Fig. 1]{Mor:09a};
        $\op{{\ell}G4_1}$ and $\op{{\ell}G4_2}$ are entries $\op{G4_4}$ and $\op{G4_5}$ in
        \cite[Table 4]{BePaPaWa:23}.
        \label{fig:graphs2plus}
        }
\end{figure}

\begin{remark}
\label{rem:obstructions}
A genus-$3+$ (irreducible) factor can have 0 crossings, 
with the
\textit{tetrahedron} $\op{T0}$
being the only possibility 
(see Fig. \ref{fig:trivialgraphs}), a
spatial graph represented by the vertices and the edges of a tetrahedron, that
can be represented by a diagram with no crossings and four trivalent vertices.
The remaining four trivial spatial graphs in Fig. \ref{fig:trivialgraphs} are excluded
because the \textit{cylinder} is itself a $2$-sum (of two theta-curves) hence it has connectivity $2$, whereas the remaining
three have connectivity $1$, 
 so that the only remaining $3$-connected trivial spatial graph is $\op{T0}$.
There are no genus-$3+$ factors with one crossing and there are 
exactly two genus-$3+$ factors with
two crossings, see Fig. \ref{fig:x2-graphs} (they are \#4, gtype \textit{cylinder}, and \#3, gtype \textit{pompomhat} in \texttt{embeddings\_2}).
Direct inspection shows that \#1 and \#2 in \texttt{embeddings\_2}, the only remaining
entries in that list, are R-equivalent to \#4
If it has $3$ crossings then an analysis of \texttt{embeddings\_3}
(Fig. \ref{fig:x3-graphs} and Table \ref{tab:three_crossings_graphs})
shows that all R-minimal
diagrams representing it require
at least $2$ obstructions in order to prevent IH-moves that show reducibility or nonminimality.
\draftMP{[Add a figure!]}
\begin{table}
  \begin{center}
    \begin{tabular}{c|c|c|c|c}
  {\bf ref} & {\bf gtype} & {\bf R-min?} & {\bf obs.} & equivs. \\
  \hline
  1 & pompom & no & \\
  2 & pompom & no & \\
  3 & pompom & no & \\
  4 & tetra  & no & \\
  5 & tetra  &    & 2 \\
  6 & cyl    & no & \\
  7 & tetra  &    & 2 & \#17 \\
  8 & tetra  &    & 2 & \#17 \\
  9 & tetra  &    & 2 & \#17 \\
 10 & tetra  &    & 2 & \#17 \\
 11 & cyl    &    & 2 & \#18 \\
 12 & pompom & no & \\
 13 & tetra  &    & 2 \\
 14 & tetra  &    & 2 & \\
 15 & tetra  &    & 2 & \#14\\
 16 & pompom & no & \\
 17 & tetra  &    & 2 \\
 18 & cyl    &    & 2 \\
    \end{tabular}
  \end{center}
    \caption{Entries in file \texttt{embeddings\_3}, three crossings.
     {\bf Ref} is the reference number in that file,
     {\bf obs.} is the number of obstructions required to make the entry either
     reducible or non-minimal.
     We note that \#4 is a link with type $[1,0,1]$, all the other entries are connected
     (type $[0,0,1]$).
     There are at most $5$ R-minimal, irreducible and inequivalent spatial graphs
     with three crossings that are not the result of a $3$-sum, see Fig. \ref{fig:x3-graphs}.}
    \label{tab:three_crossings_graphs}
\end{table}
In other words we must add at least two more factors to different arcs.
More specifically: in \#13 and \#17 we see two regions corresponding to Figure \ref{fig:moreautononminimal}(left)
with $2$ trivalent vertices involved, so that we must obstruct two IH-moves to avoid nonminimality.
In \#5 we see three partially overlapping patches corresponding to Figure \ref{fig:moreautononminimal}(left)
each involving $3$ trivalent vertices leading to nonminimality if we allow for two IH-moves for
at least one of the three patches;
we must thus obstruct an IH-move at at least two of the three internal segments.
In \#14 the patch containing the small triangle becomes a $1$-sum with two IH-moves, thus requiring two
obstruction to avoid being equivalent to a $1$-connectivity diagram, a situation covered in
Subsection \ref{sec:typeone}.
Finally in \#18 one obstruction is clearly required at the top horizontal segment;
to show that also the lower horizontal segment must be obstructed observe that
using IH moves on that segment allows to revert the overpasses at the two
lower crossings, we can then apply Reidemeister move R3 followed by R1 to reduce
the number of crossings.

At present we do not have a complete list of distinct genus-$3+$ factors with $4$ crossings.
\end{remark}

\begin{remark}
A genus-$2+$ (irreducible) factor can have 0 crossings, the only possibility being the trivial theta curve,
otherwise it has at least $2$ crossings $\op{hG2}$ ($2_1$ in \cite{Mor:07b}) or more
($\op{tG3}$: $3_1$ in \cite{Mor:09a}, $\op{hG4}$: $4_1$ in \cite{Mor:07b},
$\op{tG4}$: $4_1$ in \cite{Mor:09a}, ...), see Fig. \ref{fig:graphs2plus}.
\end{remark}

$\bullet$ 
{\sl Case $n = 5$.}

We cannot have more than $3$ classical links (otherwise the
total number of crossings would be larger than $6$). Moreover, less then $3$ classical links
would imply more then two factors with genus $> 1$, impossible.
The only possibility is then to have $2$ genus-$2+$ and $3$ links.
Each link contributes  at least $2$ crossings, 
therefore we have necessarily 2 theta-curves ($\op{tG0}$)
and three Hopf links $\op{L2a1}$.
However each theta curve requires at least three obstructions to prevent IH-moves that would
render the diagram equivalent to a $1$-connectivity diagram with no crossings increase (hence included
in $\Done$, see start of this Section).
Even if we directly sum the two
theta curves (obtaining the \textit{cylinder} of Fig. \ref{fig:trivialgraphs}) we still need $4$ obstructions, too many. \draftGB{
[cioe' per evitare riducibilita'?]}\draftMP{[Per la precisione per evitare l'equivalenza con un diagramma $1$-connesso senza
aumentare il numero di incroci]}
In conclusion we cannot have any instance with $n=5$.

$\bullet$ 
{\sl Case $n=4$.}

%\draftMP{[The following text seems not to be required:
%There are two possible tree structures representing the set of $3$ $2$-sums, one is a chain of $4$ objects, the
%other consists of a central object added to each of the other three objects.]}

{\sl Subcase $n = 4$, one $\fthree$ (a genus-$3+$ graph) and three 
$\flink$ (standard knots/links).}

Each classical link contributes for at least $2$ crossings, so that
we are forced to have $\fthree = \op{T0}$ (the \textit{tetrahedron}) and three $\op{L2a1}$ (Hopf links).
In this way we can obstruct at most 3 IH-moves for
the tetrahedron. However an IH-move at a non-obstructed edge of $\op{T0}$ allows to obtain a cylindrical
diagram with no crossings which requires 4 obstructions, too many.

{\sl Subcase $n = 4$, two $\ftwo$ and two $\flink$.}

The two $\flink$ contribute for at least 4 crossings, 
so that the two $\ftwo$ must
contribute a total of at most $2$ crossings. Hence at least one of the two $\ftwo$ is the trivial
theta curve $\op{tG0}$, which requires three obstructions, one for each other factor.
The only possibility is that the other $\ftwo$ is $\op{hG2}$ (Fig. \ref{fig:graphs2plus}), obstructing one of
the arcs of $\op{tG0}$ and being reciprocally obstructed by $\op{tG0}$ in order to avoid
a ``crossing-reducing'' IH-move.
The other two arcs of $\op{tG0}$ must be obstructed by the two 
Hopf links $\op{L2a1}$.
As a result we obtain entry $6_{11}$ in Table \ref{tab:hlinks} of type $[2,0,1]$.

$\bullet$ 
{\sl Case $n = 3$.}

{\sl Subcase $n = 3$ and two classical links.}

We have one $\fthree$ and two classical links.
The $\fthree$ cannot be the tetrahedron $T0$ (it requires too many obstructions to avoid
reducibility), hence it contributes for at least $2$ crossings.
Consequently the two classical links factors must both be an $\op{L2a1}$ forcing $\fthree$ to have
$2$ crossings, i.e. it is one of the two in Fig. \ref{fig:x2-graphs}.
Now $\fthree$ requires at least two obstructions, so that it must be directly involved with
a $2$-sum with each of the two $\op{L2a1}$.
We have two possibilities since both graphs with genus $3$ and two crossings can be obstructed
by operating on two arcs.

\begin{itemize}
\item $\fthree$ is shown in Fig. \ref{fig:x2-graphs} left.
We must obstruct the two
IH-moves on the two vertical arcs using $\op{L2a1}$, we obtain a diagram that is
equivalent to entry $6_{11}$ in Table \ref{tab:hlinks},
with a total of six crossings and type $[2,0,1]$.
\item $\fthree$ is shown in Fig. \ref{fig:x2-graphs} right.
We have two choices of which pair of arcs to obstruct.
 \begin{itemize}
 \item If we do not obstruct the IH-move involving the vertical segment in
 Fig. \ref{fig:x2-graphs} (right) we are forced to obstruct both the IH-moves involving
 the two curved external arcs.
 We end up again with a diagram equivalent to entry $6_{11}$ in Table \ref{tab:hlinks}.
 \item Otherwise we obstruct the IH-move by performing a $2$-sum with $\op{L2a1}$
 involving the vertical segment.
 At this point not obstructing the horizontal segment allows to show the existence of a
 $1$-connectivity diagram with a maximum of six crossings, that will be captured in Subsection
 \ref{sec:typeone}.
 Obstructing one of the two external curved arcs also leads to the same conclusion.
 We remain with no possibilities.
 \end{itemize}
\end{itemize}

{\sl Subcase $n = 3$ and one knot/link.}

Hence we have two $\ftwo$ factors (and one $\flink$).
We cannot have any theta-curve, since each would require three obstruction, too many
even if we take into account the possibility of obstructing with the other $\ftwo$.
Hence each $\ftwo$ and the single $\flink$ contributes at least $2$ crossings, so
that each contributes exactly $2$ crossings, otherwise we would have more than six crossings.
At this point we are forced to have two $\op{hG2}$ and one $\op{L2a1}$.
They can be combined in a chain in two possible ways
\begin{itemize}
\item $\op{hG2} \twosum \op{hG2} \twosum \op{L2a1}$:
If the two $\op{hG2}$ mutually obstruct their IH-move we first obtain entry $4_1$ in
Table \ref{tab:hlinks} and we still need to add an $\op{L2a1}$.
We have two inequivalent choices for that
\begin{itemize}
\item Performing the $2$-sum at the connecting arc of one of the handcuff components,
obtaining entry $6_8$ in Table \ref{tab:hlinks};
\item Performing the $2$-sum at a loop of one of the handcuff components,
obtaining entry $6_9$ in Table \ref{tab:hlinks}.
\end{itemize}
In both cases we end up with a handlebody link with $6$ crossings and type $[1,2]$.
Otherwise the IH-move in the central $\op{hG2}$ factor is obstructed by the
$\op{L2a1}$ and {\bf not} by the other $\op{hG2}$, forcing the placement of
$\op{L2a1}$, we end up with
\begin{itemize}
\item entry $6_5$ in Table \ref{tab:hlinks}, $6$ crossings and type $[1,0,1]$.
\end{itemize}
\item $\op{hG2} \twosum \op{L2a1} \twosum \op{hG2}$:
The central $\op{L2a1}$ can be involved with the same component for both $2$-sums, which leads to an already
covered case, entry $6_8$ in Table \ref{tab:hlinks}.
Otherwise each component of the $\op{L2a1}$ is involved in each of the two 2-sums with the effect of obtaining
\begin{itemize}
\item entry $6_2$ in Table \ref{tab:hlinks},
$6$ crossings and type $[0,2]$.
\end{itemize}
\end{itemize}

$\bullet$ 
{\sl Case $n = 2$.}

We have two possibilities, either one factor is an $\fthree$ and the other
is a knot/link, or both factors are an $\ftwo$.

{\sl Subcase $n = 2$ and no classical link factors.}
None of the two $\ftwo$ factors can be the theta-curve, since we cannot obstruct all three IH-moves.
Hence they contribute with at least $2$ crossings each.
We consider the possible total number of crossings:

\begin{itemize}
\item
{\bf $4$ crossings:}
the two factors are forced both to be $\op{hG2}$ and the $2$-sum must obstruct the IH-move
on both factors.
We obtain
$\op{hG2} \twosum \op{hG2}$, entry $4_1$ in Table \ref{tab:hlinks};
\item
{\bf $5$ crossings:}
one of the two factors is forced to be $\op{hG2}$, the other has to be $\op{tG3}$, the
arcs involved in the $2$-sum are forced in order to obstruct an unwanted IH-move, we
obtain 
$\op{tG3} \twosum \op{hG2}$, entry $5_1$ in Table \ref{tab:hknots};
\item
{\bf $6$ crossings:}
the two factors can contribute with $3$ crossings each ($\op{tG3} \twosum \op{tG3}$)
or with $2 + 4$ crossings.
In the first case ($6 = 3 + 3$) we recall that $\op{tG3}$ is chiral, so that we have the option
to take both with the same handedness, obtaining entry $6_{13}$ in Table
\ref{tab:hknots}, or with opposite handedness, obtaining entry $6_{12}$.
In the second case ($6 = 4 + 2$) we have five choices for the factor with $4$ crossings: $\op{hG4}$,
$\op{tG4}$, $\op{HK4}$, $\op{{\ell}G4_1}$ and $\op{{\ell}G4_2}$ (Fig. \ref{fig:graphs2plus}).
If not indicated otherwise the arcs involved in the $2$-sum are again forced to avoid unwanted IH-moves.
Spatial graph $\op{hG2}$ is amphichiral resulting in only one possible result for each
$2$-sum.
We consider the possible choices for the factor with $4$ crossings:
\begin{itemize}
\item
$\op{hG4} \twosum \op{hG2}$: we obtain entry $6_3$ in Table \ref{tab:hlinks};
\item
$\op{tG4} \twosum \op{hG2}$: we obtain entry $6_{10}$ in Table \ref{tab:hknots};
\item
$\op{HK4} \twosum \op{hG2}$: we have two choices for the arc of $\op{HK4}$ to use
for the $2$-sum.
If we use the arc connecting the two nodes we obtain entry $6_4$ in Table \ref{tab:hlinks}.
If we use a loop arc we obtain entry $6_{11}$ in Table \ref{tab:hknots}.
\item
$\op{{\ell}G4_1} \twosum \op{hG2}$: note that if we do not obstruct the central
horizontal segment of the theta-curve component of $\op{{\ell}G4_1}$ we can
use IH-moves and Reidemeister moves to revert the two lower crossings in
$\op{{\ell}G4_1}$ allowing to reduce the number of crossings.
The arc to use for the $2$-sum is then forced and we obtain entry $6_6$ in Table
\ref{tab:hlinks}.
\item
$\op{{\ell}G4_2} \twosum \op{hG2}$: we are forced to use bottom connecting
arc of $\op{{\ell}G4_2}$ for the $2$-sum in order to avoid a Reidemeister move
that decreases the number of crossings.
We obtain entry $6_{10}$ in Table \ref{tab:hlinks}.
\end{itemize}
\end{itemize}

{\sl Subcase $n = 2$ and one classical link.}

We have one $\fthree$ and one $\flink$.
The $\flink$ contributes for at least $2$ crossings leaving at most
$4$ crossings for $\fthree$.
If $\fthree$ has $3$ crossings or less, then we need to obstruct
at least two IH-moves at different arcs
as we already observed in Remark \ref{rem:obstructions}.

We then remain with the $\fthree$ having $4$ crossings and $\flink$
with $2$ crossings, forcing it to be $\op{L2a1}$.
To conclude we would need to build a table of R-minimal diagrams with
$4$ crossings by analizing the entries in file \texttt{embeddings\_4}
and selecting those that become minimal upon obstructing an
IH-move on a single arc.
However file \texttt{embeddings\_4} contains $110$ entries, too many
to carry out this goal.
We then opt for a different approach that we shall describe in the next
Section.
With the described construction we obtain a set of $25$ diagrams then
reduced to $18$ ($16$ with six crossings) to form the set $\DtwoA$.

We summarize the analysis of this subsection with the following lemma.

\begin{lemma}\label{lm:two_connected_IH_minimal}
An irreducible handlebody-knot admitting a
minimal diagram with connectivity $2$ and crossing number $c \leq 6$ with
the exclusion of those discussed in the next Subsection \ref{sec:typetwospecial}
is equivalent, up to mirror image, to one of the following handlebody knots: 
\begin{equation}\label{list:handlebody_knot_type_2}
 5_1, 6_{10}, 6_{11}, 6_{12}, 6_{13}.
\end{equation}
in Table \ref{tab:hknots}.

An irreducible non-split multi-component handlebody-link admitting a
minimal diagram with connectivity $2$ and crossing number $c \leq 6$ 
is equivalent, up to mirror image, to one of the following handlebody links: 
\begin{equation}\label{list:handlebody_link_type_2}
 4_1, 6_2, 6_3, 6_4, 6_5, 6_6, 6_8, 6_9, 6_{10}, 6_{11}.
\end{equation}
in Table \ref{tab:hlinks}.

\end{lemma}

%
% ====================================================================
%

\subsection{Minimal diagrams with connectivity $2$, special case}
\label{sec:typetwospecial}

The number of planar embeddings with four crossings (110) is too large for a manual analysis as was done for
up to three crossings.

However the only possible missing $2$-sums are the result of adding a $3$-connected $4$-crossings spatial
graph $G$ to a Hopf link ($\op{L2a1}$).
This sum can be done in different ways, precisely there is one $2$-sum for each
arc connecting pairs of $3$-valent nodes of the spatial graph, six choices, and one $2$-sum for each
genus-$1$ connected component.
There are $7$ out of $110$ planar embeddings with a genus-$1$ component, all of which corresponding to
a handlebody of type [1,0,1].
Direct inspection of the result of a $2$-sum of the genus-$1$ component to the Hopf link L2a1 shows that
none of the diagrams obtainable in this manner is minimal.
We then remain with the $6$ arcs connecting pair of vertices (or a loop attached to a vertex),
hence $6$ possible $2$-sums for each graph $G$.

We can then build in an automated way all such $2$-sums for all entries in file \texttt{embeddings\_4}.
Then we can build all possible choices of the overcrossing information and populate the bigtree as
done for the $3$-connected case described in Subsection \ref{sec:typethree}.

The result of the dynamic process described in Section 
\ref{sec:typethree} limited to the
diagrams constructed as described above results in the following Lemma.

\begin{lemma}\label{lm:two_connected_IH_minimal_special}
No new irreducible handlebody-knot admitting a minimal diagram with connectivity $2$ and crossing number
at most $6$ and not listed in \eqref{list:handlebody_knot_type_2}, up to mirror image, is uncovered by the
procedure illustrated in this subsection.
%
%An irreducible handlebody knot admitting a
%minimal diagram with connectivity $2$ and crossing number $6$ 
%is equivalent, if not already listed in \eqref{list:handlebody_knot_type_2}
%and up to mirror image, to one of the following handlebody knots: 
%\begin{equation}\label{list:handlebody_knot_type_2_special}
% NONE!
%\end{equation}
%in Table \ref{tab:hknots}.
%

An irreducible non-split multi-component handlebody link admitting a
minimal diagram with connectivity $2$ and crossing number $6$ 
is equivalent, if not already listed in \eqref{list:handlebody_link_type_2}
and up to mirror image, to the following handlebody-link: 
\begin{equation}\label{list:handlebody_link_type_2_special}
 6_7.
\end{equation}
in Table \ref{tab:hlinks}.

\end{lemma}

% END of \nada

Consequently, we obtain a refined set $\DtwoB$ of %
$101$ diagrams.
As for the $3$-connected case
we arrive at this final number in two steps, first all situations depicted in Figures \ref{fig:autononminimal}
and \ref{fig:moreautononminimal} with the exception of the IH+V move are recognized by the \texttt{appcontour} software
and used to generate a set of $183$ diagrams.

This set is then reduced to $101$ entries by automatically recognizing situations matching the move ``IH+V'' of Figure
\ref{fig:moreautononminimal} and, for any set of diagrams that are directly connected by IH-moves only one is retained.

\draftMP{[A few numbers:
\begin{itemize}
\item
$660 = 110 \times 6$ entries in \texttt{embeddings\_6extra}
\item
$1487$ diagrams in \texttt{pairs6extra-two\_sum}
\item
$1178$ are recognized as nonminimal by \texttt{search334.sh}
leaving $309$ diagrams
\item
$125$ are further recognized as nonminimal by \texttt{search3334.sh}
leaving $184$ diagrams
\item
$1$ recognized by \texttt{search33334.sh}, we remain with $183$
\item
$82$ are further recognized as equivalent up-to-ih to a previous one by
\texttt{search\_uptoih.sh} leaving us with $101$ diagrams
\end{itemize}
}

We finally define $\Dtwo = \DtwoA \cup \DtwoB$ with all the constructed diagrams having arc connectivity $2$.

%
% ====================================================================
%

\subsection{Minimal diagrams with connectivity 1}
\label{sec:typeone}

See Table \ref{tab:typeone} and Fig. \ref{fig:typeonedetail}.

\subsubsection{Double $1$-sum (doubly reducible handlebody-links).}
\label{sec:doubletypeone}

Each $1$-sum increases the number of trivalent nodes in the diagram by
$2$ (i.e. it decreases the Euler characteristic by $1$).
It is possible then to have two $1$-sums between standard knots/links, so
we need to find all sets of three classical links with total number of crossings
at most $6$, also including non-prime classical links and taking into account
the handedness of the factors.

Choosing the three ingredients does not give a single (or two if there are
two chiral factors) result when we deal with the nonconnected case.  Even the
order with which we perform the operations becomes important since the central
factor can be involved with the other two with distinct connected components.

We have only two nonprime knots with at most $6$ crossings, namely the
square knot (the $2$-sum of two trefoil knots K3 with opposite handedness)
and the granny-knot ($2$-sum of two K3 with the same handedness).  They
must be taken into account as a possible factor to be $1$-added with two
unknots (K0) leading to two distinct results.
However the nonconnected case is much more involved:
L2\#L2 (four crossings) is a chain of three rings, but L2\#L2\#L2 gives two
distinct composite links, one is a chain of four rings, the other consists of
a central ring interlaced with three other rings.
The result must then be $1$-added with two unknots which requires a choice of
two connected components leading to a total of $10$ distinct configurations
of L2\#L2\#L2--K0--K0.
The total number of doubly reducible handlebody-links is $73$.

% ============================================
\begin{table}
  \begin{center}
  \begin{tabular}{c|ccccccc|}
  \textbf{crossings} & [0,0,1] & [0,2] & [1,0,1] & [1,2] & [2,0,1] & [2,2] & [3,0,1] \\
  \hline
  0 & 1(1) &      &      &      &      &       &      \\
  2 &      & 1(1) & 1(1) &      &      &       &      \\
  3 & 1(1) &      &      &      &      &       &      \\
  4 & 2(1) & 2(1) & 2(1) & 3(3) & 3(3) &       &      \\
  5 & 6(2) & 4(3) & 5(4) &      &      &       &      \\
  6 &23(7) &13(6) &15(6) &19(8) &12(5) &13(12) & 8(7) \\
  \hline
  \textbf{total} & 
     33(12) &20(11) &23(12) &22(11) &15(8) &13(12) & 8(7) 
  \end{tabular}
  \end{center}
  \caption{Number of reducible handlebody-links up to six crossings by their type.
           Inside parentheses the number of doubly reducible objects.}
  \label{tab:typeone}
\end{table}

\subsubsection{Single $1$-sum.}
\label{sec:singletypeone}

An analysis similar to that of the previous Subsection \ref{sec:doubletypeone} can be
performed, with the novelty now that one of the two factors must have Euler characteristic
$-1$, i.e. genus $1$ plus number of connected components.
We obtain a total of $61$ distinct results for a grand total of $134$.

Collecting all diagrams of Subsections \ref{sec:doubletypeone} and \ref{sec:singletypeone} we obtain a set $\Done$ of
(potentially distinct) $134$ diagrams.

Whether the elements of set $\Done$ are indeed mutually distinct is not guaranteed.
We shall treat them together with the other sets $\Dthree$ and $\Dtwo$, using the dynamic process
discussed in Section \ref{sec:dynamicapproach} and try to distinguish them by computing suitable
invariants.
Indeed we find a hard-pair (i.e. a pair of in principle distinct reducible handlebody-links that we
could not distinguish using computed invariants).
This hard-pair will be discussed in Section \ref{sec:hardpairs}.

\begin{figure}
\begin{center}
\begin{subfigure}{.32\linewidth}
\begin{tabular}{|c|cc|}
\multicolumn{3}{c}{\textbf{[0,0,1]}}\\
\hline
\textbf{x} & \textbf{n} & \textbf{name} \\
\hline
0 & 1 & K0--K0--K0 \\
3 & 1 & K3--K0--K0 \\
4 & 1 & HK4--K0 \\
4 & 1 & K4--K0--K0 \\
5 & 2 & K5$_{[1-2]}$--K0--K0 \\
5 & 4 & HK5$_{[1-4]}$--K0 \\
6 & 2 & K3--K3--K0 \\
6 & 2 & K3\#K3--K0--K0 \\
6 & 3 & K6$_{[1-3]}$--K0--K0 \\
6 & 16 & HK6$_{[1-16]}$--K0 \\
\hline
\multicolumn{3}{c}{ }\\
\multicolumn{3}{c}{ }\\
\end{tabular}
\end{subfigure}
\begin{subfigure}{.32\linewidth}
\begin{tabular}{|c|cc|}
\multicolumn{3}{c}{\textbf{[0,2]}}\\
\hline
\textbf{x} & \textbf{n} & \textbf{name} \\
\hline
2 & 1 & L2--K0--K0 \\
4 & 1 & HL4--K0 \\
4 & 1 & L4--K0--K0 \\
5 & 1 & HL5--K0 \\
5 & 1 & K3--L2--K0 \\
5 & 1 & K3\#L2--K0--K0 \\
5 & 1 & L5--K0--K0 \\
6 & 1 & K4--L2--K0 \\
6 & 1 & K4\#L2--K0--K0 \\
6 & 3 & L6a[1-3]--K0--K0 \\
6 & 8 & HL6$_{[1-8]}$--K0 \\
\hline
\multicolumn{3}{c}{ }\\
\end{tabular}
\end{subfigure}
\begin{subfigure}{.32\linewidth}
\begin{tabular}{|c|cc|}
\multicolumn{3}{c}{\textbf{[1,0,1]}}\\
\hline
\textbf{x} & \textbf{n} & \textbf{name} \\
\hline
2 & 1 & L2--K0--K0 \\
4 & 1 & HL4--K0 \\
4 & 1 & L4--K0--K0 \\
5 & 1 & HL5--K0 \\
5 & 1 & K3--L2--K0 \\
5 & 2 & K3\#L2--K0--K0 \\
5 & 1 & L5--K0--K0 \\
6 & 1 & HK4--L2 \\
6 & 1 & K4--L2--K0 \\
6 & 2 & K4\#L2--K0--K0 \\
6 & 3 & L6a[1-3]--K0--K0 \\
6 & 8 & HL6$_{[1-8]}$--K0 \\
\hline
\end{tabular}
\end{subfigure}
\end{center}
\begin{center}
\begin{subfigure}{.32\linewidth}
\begin{tabular}{|c|cc|}
\multicolumn{3}{c}{\textbf{[1,2]}}\\
\hline
\textbf{x} & \textbf{n} & \textbf{name} \\
\hline
4 & 1 & L2--L2--K0 \\
4 & 2 & L2\#L2--K0--K0 \\
6 & 1 & HL4--L2 \\
6 & 2 & L4--L2--K0 \\
6 & 3 & L4\#L2--K0--K0 \\
6 & 2 & L6a[4-5]--K0--K0 \\
6 & 1 & L6n1--K0--K0 \\
6 & 1 & HL6$_{9}$--K0 \\
6 & 2 & HL6$_{10}$--K0 \\
6 & 2 & HL6$_{11}$--K0 \\
6 & 2 & HL6$_{12}$--K0 \\
6 & 1 & HL6$_{13}$--K0 \\
6 & 2 & HL6$_{14}$--K0 \\
\hline
\end{tabular}
\end{subfigure}
\begin{subfigure}{.32\linewidth}
\begin{tabular}{|c|cc|}
\multicolumn{3}{c}{\textbf{[2,0,1]}}\\
\hline
\textbf{x} & \textbf{n} & \textbf{name} \\
\hline
4 & 1 & L2--L2--K0 \\
4 & 2 & L2\#L2--K0--K0 \\
6 & 1 & HL4--L2 \\
6 & 1 & L4--L2--K0 \\
6 & 3 & L4\#L2--K0--K0 \\
6 & 1 & L6n1--K0--K0 \\
6 & 6 & HL6$_{[9-14]}$--K0 \\
\hline
\end{tabular}
\end{subfigure}
\begin{subfigure}{.32\linewidth}
\begin{tabular}{|c|cc|}
\multicolumn{3}{c}{\textbf{[2,2]~~x=6}}\\
\hline
\textbf{n} & & \textbf{name} \\
\hline
 1 & \multicolumn{2}{|c|}{L2--L2--L2}\\
 5 & \multicolumn{2}{|c|}{L2\#L2--L2--K0}\\
 6 & \multicolumn{2}{|c|}{L2\#L2\#L2--K0--K0}\\
 1 & \multicolumn{2}{|c|}{HL6$_{15}$--K0}\\
\hline
\multicolumn{3}{c}{ }\\
\multicolumn{3}{c}{\textbf{[3,0,1]~~x=6}}\\
\hline
\textbf{n} & & \textbf{name} \\
\hline
 1 && \!\!\!\!L2--L2--L2 \\
 2 && \!\!\!\!L2\#L2--L2--K0 \\
 4 & \multicolumn{2}{|c|}{L2\#L2\#L2--K0--K0}\\
 1 && \!\!\!\!HL6$_{15}$--K0 \\
\hline
\end{tabular}
\end{subfigure}

\end{center}
\caption{List of reducible non-split handlebody-links
with up to $6$ crossings for the relevant values of type.
{\bf x} is the crossing number, {\bf n} is the number of possible
distinct results for the expression in the last column.
Abbreviations used: K$n_i$ and L$n*_i$ refer to an entry in the Rolfsen knot table and the
Thistlethwaite link table respectively; HK$n_i$, HL$n_i$ are handlebody-knots/links in
\cite{IshKisMorSuz:12} and \cite{BePaPaWa:23}.
The subscript is omitted whenever there is no confusion possible.
The notation ``[i-j]'' is a range of integers and the corresponding object should be
considered multiple times.
Knot K3 used twice in the same expression produces two results based on the choice of
their handedness.}
\label{fig:typeonedetail}
\end{figure}

%
% ===============================================================================
%

\subsection{Conclusion}\label{sec:completenessconclusion}
Let $\Dcal = \Done \cup \Dtwo \cup \Dthree$.

\begin{theorem}[Completeness]\label{teo:completeness}
For every non-split handlebody-link $V$ with $\chi(V) = -2$ with crossing number six or less, there is a minimal diagram in $\Dcal$ representing $V$. 
\end{theorem}
\begin{proof}
As shown in Section \ref{sec:typeone}, if $V$ has a minimal diagram with connectivity $1$, then $\Done$ contains a minimal diagram of $V$.        
Suppose $V$ has no minimal diagram with connectivity $1$. If $V$ has a minimal diagram with connectivity $2$, Sections \ref{sec:typetwo} and
\ref{sec:typetwospecial} imply $\Dtwo$ contains a minimal diagram of $V$; otherwise
$V$ has a minimal diagram in $\Dthree$; see Section \ref{sec:typethree}.
\end{proof}

\begin{remark}
In general, it is not known whether the minimal diagram of a $1$-sum always has connectivity $1$.
Similarly, it is unclear whether a minimal diagram of a handlebody-link having a spine that is a $2$-sum always has connectivity $2$.

\draftMP{[Should we say something about split handlebodies?]}

The enumeration is, however, based on the connectivity of minimal diagrams, and not on the decomposition type of handlebody-knots, and therefore does not depend on the aforementioned open problems. For instance, if a $1$-sum has no minimal diagrams with connectivity $1$, then its minimal diagram(s) can be found in $\Dtwo \cup \Dthree$.
\end{remark}

%%%%%%%%%%%%%%%%%%%%%%%%%%%%%%%%%%%%%%%%%%%%%%%%%%%%%%%%%%%%%%%%%%%%%%%%%%%%%

%%%%%%%%%%%%%%%%%%%%%%%%%%%%%%%%%%%%%%%%%%%%%%%%%%%%%%%%%%%%%%%%%%%%%%%%%%%%%
\section{No duplicates in the tables}\label{sec:uniqueness}
Our goal is to partition
 the diagrams in $\Dcal$ into classes, 
each containing only diagrams
that represent the same handlebody-knot, up to equivalence and mirror image.
Given the large size of $\Dcal$ we adopt a dynamic strategy in order to achieve this goal.

\subsection{A dynamic approach}
\label{sec:dynamicapproach}

The so-called `bigtree' is a dynamic data structure, actually a {\it growing} mathematical rooted tree
with data attached to each vertex.
It initially consists of a degenerate tree with a single vertex that is also a leaf and contains
diagrams of handlebody-links up to a given number of crossings.
It is important that it contains at least one minimal diagram for 
all handlebody-links (up
to mirror image), including all reducible unsplit handlebodies.

\draftMP{[Probably it is better not to mention split diagrams in the bigtree, removed a few lines]}
\nada{
Although not strictly necessary it might be useful to also include
a few split handlebodies in order to recognize connected diagrams that represent a split
handlebody.
It seems reasonable to bound the number of crossings for split diagrams to four.
Actually it turned out that none of the diagrams in $\Dcal$ represents a split handlebody.
}

Tables \ref{tab:hknots} and \ref{tab:hlinks} list \textit{irreducible} handlebodies, however
the dynamic process also produces all reducible unsplit handlebody-links.
% and some split handlebody-links.
Populating the bigtree with also reducible objects is
key in order to recognize possible
% {\it splittable} or
 reducible handlebodies.
\draftMP{[Al momento non mi ricordo più il motivo di includere oggetti "splitted" (fino a 4 incroci)]}

The goal of the growing process of the tree is to 
reach a final state where all leaves
either contain all equivalent diagrams, or there are a small number of, usually $2$, distinct
equivalence classes that are difficult to distinguish by computable invariants but are nevertheless
proven to be inequivalent.
An example could be a knot sum of two chiral handlebody-link.
A snapshot of the contents of the bigtree is available online at \cite{bigtable:web}.

At an intermediate state each non-leaf vertex of the bigtree is tagged with an {\it invariant} and
all (arcs) branches 
directed away from the root correspond to different achieved values of that invariant.
There are typically at least two such branches.

Leaves contain a subset of the original population of diagrams 
and can be in a non-final state if we
are still unable to subdivide the 
diagrams they contain into equivalence classes (typically one, rarely two).
Otherwise, such a leaf is tagged as final.

At each step, we select a non-final leaf vertex 
and seek a computable invariant that is capable of discriminating
among the diagrams it contains, while 
ensuring that computing the invariant for all diagrams in
the (former) leaf is not prohibitively time-consuming.
We then tag such vertex with the selected invariant, add to it new departing branches tagged with the
resulting invariant values and move all diagrams in the former leaf to the relevant newly created leaf.

Also, we can also always drop diagrams that are recognized to be nonminimal or sometimes diagrams that an automated tool can show to
be equivalent to others in the same leaf vertex.
This often happens when two diagrams differ only in the
over/under crossings information.

The reason for the whole construction is that 
we are dealing with a large number of
diagrams in the first place, and computing strong invariants that can sufficiently discriminate is
computationally not affordable if applied to the whole original set.

%
% ==============================================================
%

\subsection{Invariants}\label{sec:invariants}

During the dynamic process described above we used a number of invariants that we are able to compute
using the \texttt{appcontour} software \cite{appcontour}.
A substantial number of invariants, computed on links, are actually a combination of
other invariants computed on various subsets of connected components of the handlebody-link.
Here is a comprehensive list.

\begin{itemize}
\item Invariants used to distinguish between irreducible handlebody-links:
\begin{itemize}
\item \texttt{type}:
Computes the \textit{type} of the handlebody associated to the diagram.
This is the first computed invariant on the whole set $\Dcal$ 
in \eqref{eq:D} and allows e.g.
to discriminate between handlebody-knots and handlebody-links.
\item \texttt{ks\_$G$}:
The Kitano-Suzuki invariant.
Given a handlebody-link $V$, let $G$ be a finite group; 
the Kitano-Suzuki invariant 
\text{ks\_$G$} 
is
the number of conjugate classes of homomorphisms from the fundamental 
group of the complement $E(V)$ of $V$ to
$G$; see \cite{KitSuz:12}.
The choices of $G$ that we used are $\op{SL}_{2,p}$, $p$ a small prime (we used $p\in\{3,5\}$),
$G = S_n$, the symmetric group with $n \leq 6$ and $G = 
A_n$, the alternating group, with $n \in \{4,5\}$;
\item \texttt{$G$-image}: For a given handlebody-link $V$ with complement $\Compl V$, this invariant, introduced in \cite{BePaWa:20a}, is constructed as follows:
 consider the homomorphisms
 \[  
 {\iota_1}_* :
        \pi_1(\partial V)\rightarrow \pi_1(V),\quad
        {\iota_2}_* :
        \pi_1(\partial V)\rightarrow \pi_1(\Compl V)
 \]
 induced by the inclusions $\iota_1 : \partial V \to V$ and $\iota_2 : \partial V \to \Compl V$.
 Then a surjective homomorphism from $\pi_1(\Compl V)$ to a finite group $G$ is said to be \textit{proper} if it
 becomes non-surjective after precomposed with ${i_2}_*$. The $G$-image is then given by the set of the images of
 the subgroup ${i_2}_\ast(\operatorname{Ker}({i_1}_\ast))$ in $G$ under proper homomorphisms.
 We used this invariant with $G = S_6$ to distinguish between entries $6_1$ and $6_4$ of Table \ref{tab:hknots};
 \footnote{In both cases the number of generic ($572297$), surjective ($459418$) and proper (only $982$)
% the software says $113861$ but counts also the nonsurjective
 homomorphisms is the same for the
 two entries, however in $16$ cases the normal closure of the image has order $72$ in one leaf and $36$ in the other.%
%< Proper homomorphism: surface group: 72 - normal closure: 72 ("(S3 x S3) : C2")
%> Proper homomorphism: surface group: 72 - normal closure: 36 ("S3 x S3")
% Relevant files can be found in directory other_material in this Svn project
 }
\end{itemize}
\item Invariants used to distinguish between reducible handlebody-links.
Let $V$ be a handlebody-link of type $[n_1, n_2, n_3]$ with components $V_1,\dots, V_m$, $m = n_1+n_2+n_3$
(and $n_2 + 2 n_3 = 2$):
\begin{itemize}
\item $\mathcal S$-$i$: when $i$ is one of the above mentioned invariants and $\mathcal S$ is a
 suitable subset of the subsets of the connected components of $V$, then this invariant computes
 $$\sum_{\sigma \in \mathcal{S}} i \left( \bigcup_{W \in \sigma} W \right)
 $$
 (the sum might be substituted by any other commutative operation) and $\mathcal S$, depending on the type of $V$,
 must be closed under any permutation $\pi$ of the connected components such that $g(\pi(W)) = g(W)$.
 A few examples on how to chose $\mathcal S$ might be:
 \begin{itemize}
 \item $\mathcal S$ is the singleton of the set of connected components of $V$ having genus $>1$.
  That is, the invariant is computed on $V$ after removal of all 
classical knot components;
  \draftMP{[This is \texttt{H-} in the bigtree]}
 \item $\mathcal S$ contains the complement of the above, i.e. is the singleton of the set of connected components of $V$ having genus $=1$;
  \draftMP{[This is \texttt{L-} in the bigtree]}
 \item $\mathcal S$ is the set of singletons each one of those containing one of the connected components of $V$.
  \draftMP{[This is \texttt{HL-} in the bigtree]} % [1,0,1], [2,0,1], [3,0,1]
 \item $\mathcal S$ is the set of singletons each one of those containing a connected component of genus $g \geq 2$;
  \draftMP{[This is \texttt{HS-} in the bigtree]}
 \item $\mathcal S$ is the set of subsets obtained by taking one component of genus $g > 1$ and all components
  of genus $g = 1$.
  This choice allows e.g. to distinguish the 
incarnation of L2a1--L2a1--K0 having type [1,2]
  from that of HL6$_9$--K0 having the same type, where HL6$_9$ is entry $6_9$ of
  \cite[Table 1]{BePaPaWa:23};
  \draftMP{[This is \texttt{HK-} in the bigtree]}
 \item $\mathcal S$ is the set of subsets obtained by taking all components of genus $g > 1$ and all but one components of
  of genus $g = 1$;
  \draftMP{[This is \texttt{HH-} in the bigtree]} % [1,2], [2,2]
 \item $\mathcal S$ is the set of subsets obtained by removing one connected component regardless of its genus;
  \draftMP{[This is \texttt{HLX-} in the bigtree]} % [1,0,1], [2,0,1], [3,0,1]
 \end{itemize}
\item \texttt{linkingnumber}: \draftMP{[\texttt{L-linkingnumber} in the bigtree]} % [1,0,1], [2,0,1], [3,0,1]
 After removal of all genus $g>1$ components we remain with a 
classical link $L$.
 This invariant then computes and adds together the absolute value of all linking numbers between pairs of components of $L$;
\item \texttt{longitude-}\textit{i}: \draftMP{[actually: L-longitude]}
 Let $\pi$ be
the fundamental group of the handlebody-link at hand $V$ (actually a presentation of it).
 For each component $K$ of genus $g = 1$ (a classical knot) we add a preferred longitude of $K$ as a relator to $\pi$
 and compute invariant $i$ (typically the Kitano-Suzuki 
\text{ks\_$G$} 
 for some group $G$) on the resulting group.
 We then add all results together for the final value of the invariant.
 This invariant allows us e.g. to distinguish among various incarnations of L4a1\#L2a1--K0--K0
 having type [2,0,1].
\end{itemize}
\end{itemize}

\subsection{Outcome}\label{sec:outcome}
At the end of the dynamic process we arrive at a set of \textit{final leaves} of the bigtree where
no aforementioned computable invariants can further distinguish among entries.
We then manually checked, with two exceptions, using the generalized 
Reidemeister moves I to V and IH-moves (see Figs.
\ref{fig:reidemeister_a}, \ref{fig:reidemeister_b}, \ref{fig:IH-move})
 that diagrams in each final leaf are all mutually equivalent, up to mirror image.
For each of the two exceptional leaves, we partitioned the diagrams therein in two subsets
of mutually equivalent diagrams. 
The handlebody-links represented by these two subsets are called \textit{hard-pairs}.
They cannot be distinguished by any of the listed computational invariants, yet we were not able to show
their equivalence.
The two hard pairs are:
\begin{itemize}
\item $6_{12}$ and $6_{13}$ in Table \ref{tab:hknots},
\item The $1$-sums of HL6$_{10}$ and K0 given by the two genus one components of HL$6_{10}$, where HL6$_{10}$ is $6_{10}$ in
\cite[Table 1]{BePaPaWa:23}; see Fig.\ \ref{fig:typeonedetail}.
\end{itemize}
Their inequivalence within each hard pair is proved with other means in Section \ref{sec:hardpairs}. The end result is Tables \ref{tab:hknots},\ref{tab:hlinks}, and Figure \ref{fig:typeonedetail}.

%%%%%%%%%%%%%%%%%%%%%%%%%%%%%%%%%%%%%%%%%%%%%%%%%%%%%%%%%%%%%%%%%%%%%%%%%%%%%

%%%%%%%%%%%%%%%%%%%%%%%%%%%%%%%%%%%%%%%%%%%%%%%%%%%%%%%%%%%%%%%%%%%%%%%%%%%%%%
\section{Irreducibility and non-splittability}\label{sec:irreducibility}

\subsubsection*{Non-splittability}
\begin{figure}[h]
\includegraphics[scale=.1]{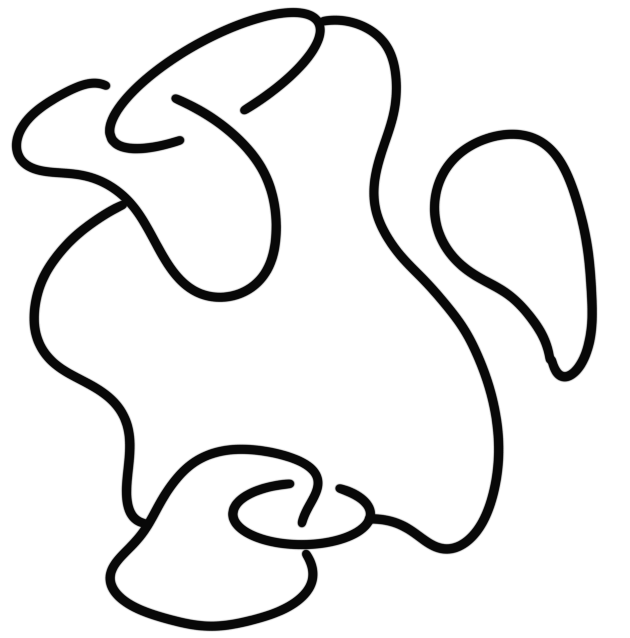}
\quad
\includegraphics[scale=.1]{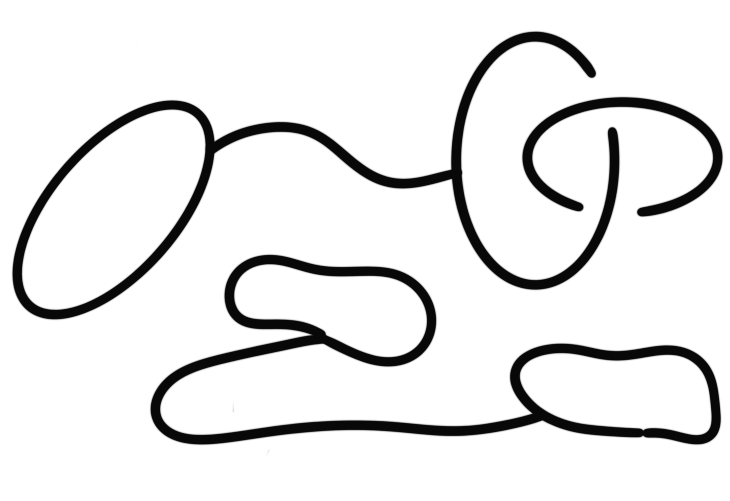}
\quad
\includegraphics[scale=.1]{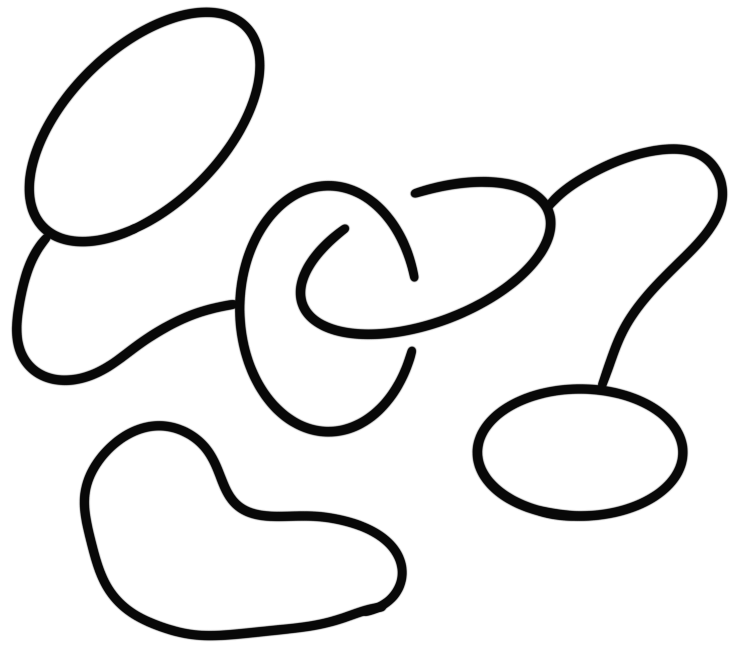}
\qquad
\includegraphics[scale=.1]{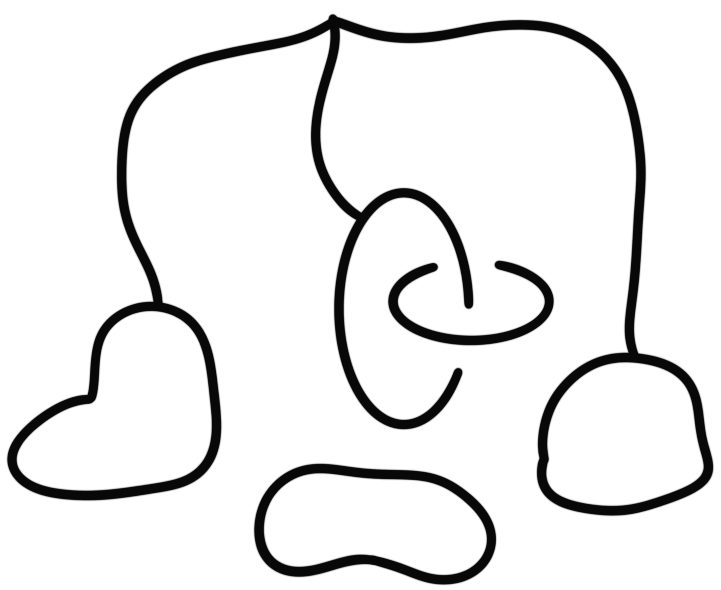}
\caption{Split handlebody-links.}
\label{fig:fake_sixeight_to_eleven}
\end{figure}

If any of $4_1$ to $6_7$ in Table \ref{tab:hlinks} is split, then 
it is equivalent to a trivial handlebody-link, and hence the fundamental group of its exterior is free and hte $\texttt{ks}_{A_5}$-invariant is $216341$, contradicting Table \ref{tab:reducibility02_101and12_201} (left).
On the other hand, if any of $6_8$ to $6_{11}$ in Table \ref{tab:hlinks} is split, then it is either trivial or equivalent to one of the four split handlebody-links in Fig.\ \ref{fig:fake_sixeight_to_eleven}. In particular, its fundamental group is either free or  
\[
\langle a,b,c,d,e  \mid aba^{-1}b^{-1}=1\rangle,
\ \text{or}\ 
\langle a,b,c,d,e  \mid\ aba^{-1}b^{-1}cdc^{-1}d^{-1}=1 \rangle,  
\]
so its $\texttt{ks}_{A_4}$-invariant is $21154$, $7330$, or $5794$, respectively,
contradicting Table \ref{tab:reducibility02_101and12_201} (right).

\subsubsection*{Irreducibility}
Recall from \cite{BePaWa:20} that
 a handlebody-link $V$ is \emph{$r$-generator} if the fundamental group of its exterior is of \emph{rank $r$}, namely, the minimal number of generators needed. Also, the \emph{genus} of a handlebody-link is defined 
as the sum $g$ of the genera of its components; in particular, a handlebody-link of type $[n_1,\dots, n_k]$ is of genus $\sum_{i=1}^k in_i$. Set $k_5:=\texttt{ks}_{A_5} (V)$ and $k_4:=\texttt{ks}_{A_4} (V)$. Then it follows from \cite[Theorems $1.1$-$2$, Corollary $1.3$, Table $1$]{BePaWa:20} the irreducibility criteria:

\begin{theorem}\label{teo:irre_criteria}
Let $V$ be an $r$-generator handlebody-link $V$ of genus $g$.   
\begin{enumerate}[label=(\roman*)]
\item\label{itm:irre_001} If $V$ has one component, $r=g+1$, and 
$k_5 \not\equiv 15 \bmod 60$, then $V$ is
irreducible.
\item\label{itm:irre_02_101} If $V$ is of type $[0,2],[1,0,1]$, $r=g$, and $k_5 \not\equiv 41 \bmod 60$, 
then $V$ is
irreducible.
\item\label{itm:irre_12_201} If $V$ is of type $[1,2],[2,0,1]$, $r=g$, and $k_4 \not\equiv 10 \bmod 12$, 
then $V$ is irreducible.
\end{enumerate}
\end{theorem} 

\begin{corollary}
All entries in Tables \ref{tab:hknots}, \ref{tab:hlinks} are irreducible.
\end{corollary}
\begin{proof}
All entries in Table \ref{tab:hknots} are of type $[0,0,1]$ and $4$-generator of genus $3$, so Theorem \ref{teo:irre_criteria}\ref{itm:irre_001} applies.
Table \ref{tab:reducibility001} lists their $k_5$'s and show their irreducibility.

In Table \ref{tab:hlinks}, $4_1, 6_1,\dots,6_7$ are of type $[0,2]$ or $[1,0,1]$ with $g=r=4$. Thus Theorem \ref{teo:irre_criteria}\ref{itm:irre_02_101} applies and their irreducibility follows from Table \ref{tab:reducibility02_101and12_201} (left). 

The remaining entries, namely $6_8, 6_9, 6_{10}, 6_{11}$ in Table \ref{tab:hlinks}, are of type $[1,2]$ or $[2,0,1]$, and all have $g=r=5$. Their irreducibility follows from Table \ref{tab:reducibility02_101and12_201} (right) by Theorem \ref{teo:irre_criteria}\ref{itm:irre_12_201}.   
\end{proof}

\begin{table}
\begin{tabular}{ccr|ccr}
entry & $k_5$ & $\bmod 60$ & entry & $k_5$ & $\bmod 60$ \\
\hline
$5_1$  & 4988 &  8 & $6_7$  & 4111 & 31 \\
$6_1$  & 4273 & 13 & $6_8$  & 4988 &  8 \\
$6_2$  & 4111 & 31 & $6_9$   & 4273 & 13\\
$6_3$  & 4273 & 13 & $6_{10}$  & 4553 & 53 \\
$6_4$  & 4273 & 13 & $6_{11}$ & 8313 & 33\\
$6_5$  & 7087 &  7 & $6_{12}$ & 5276 & 56 \\
$6_6$  & 6956 & 56 & $6_{13}$ & 5276 & 56 \\
\end{tabular}
\caption{Type $\text{[0,0,1]}$; 
the third column shows the result of $k_5 \bmod 60$.
% that should be compared with the value $15$.
%The first column refers to Table \ref{tab:hknots}.
}\label{tab:reducibility001}
\end{table}

% ==================================================================================

\begin{table}
\begin{tabular}{ccr|ccr}
entry & $k_5$ & $\bmod 60$ & entry & $k_4$ & $\bmod 12$ \\
\hline
$4_1$ &  5110 & 10 & $6_8$     & 3746 & 2 \\
$6_1$ &  9050 & 50 & $6_9$     & 2306 & 2 \\
$6_2$ & 45015 & 15 & $6_{10}$  & 2306 & 2 \\
$6_3$ &  8474 & 14 & $6_{11}$  & 3074 & 2 \\
$6_4$ & 23170 & 10 \\
$6_5$ & 44336 & 56 \\
$6_6$ & 22136 & 56 \\
$6_7$ & 22856 & 56 \\
\end{tabular}
\caption{Left: type %$\text{htype}(V) \in 
$\{\text{[0,2]},\text{[1,0,1]}\}$;  
%The first column refers to Table \ref{tab:hlinks}, 
the third column shows the result of $k_5 \bmod 60$. 
%that should be compared with the value $41$.%
Right: type $\{\text{[1,2]},\text{[2,0,1]}\}$; %$v = \texttt{ks\_A4} (V)$.
%The first column refers to Table \ref{tab:hlinks},
the third column shows the result of $k_4 \bmod 12$.
% that should be compared with the value $10$.
}\label{tab:reducibility02_101and12_201}
\end{table}

%
%\begin{table}
%\begin{tabular}{ccr}
%entry & $k_4$ & $\bmod 12$ \\
%\hline
%$6_8$     & 3746 & 2 \\
%$6_9$     & 2306 & 2 \\
%$6_{10}$  & 2306 & 2 \\
%$6_{11}$  & 3074 & 2 \\
%\end{tabular}
%\caption{Type $\{\text{[1,2]},\text{[2,0,1]}\}$; %$v = \texttt{ks\_A4} (V)$.
%%The first column refers to Table \ref{tab:hlinks}, 
%the third column shows the result of $k_4 \bmod 12$.
%% that should be compared with the value $10$.
%}
%\label{tab:reducibility12and201}
%\end{table}
%

%%%%%%%%%%%%%%%%%%%%%%%%%%%%%%%%%
\nada{
Excerpt from a mail of Yi-Sheng about the irreducibility criteria:

\begin{verbatim}
By the way, I checked the irreducibility test for the other types; it goes as
follows:

[...]

2. [1,2], [2,0,1] belong to the same category,

if r=5 and ks_A4 is not 10 (mod 12), then it is irreducible.

In the second case, if r=5, then it seems we are done, while in the first
case, for HHL6_2 and HHL6_3, we might need ks_A5, yet even if it fails to have
a 4-generator presentation, there is another way to see its irreducibility
(Theorem 4.1 in the irreducibility test paper).
\end{verbatim}
}
%%%%%%%%%%%%%%%%%%%%%%%%%%%%%%%%%

%%%%%%%%%%%%%%%%%%%%%%%%%%%%%%%%%%%%%%%%%%%%%%%%%%%%%%%%%%%%%%%%%%%%%%%%%%%%%%

%%%%%%%%%%%%%%%%%%%%%%%%%%%%%%%%%%%%%%%%%%%%%%%%%%%%%%%%%%%%%%%%%%%%%%%%%%%%%
\section{Hard pairs}\label{sec:hardpairs}
\subsection{Irreducible hard pair}
\begin{figure}[t]
\begin{subfigure}{.3\linewidth}
\centering 
\includegraphics[scale=.1]{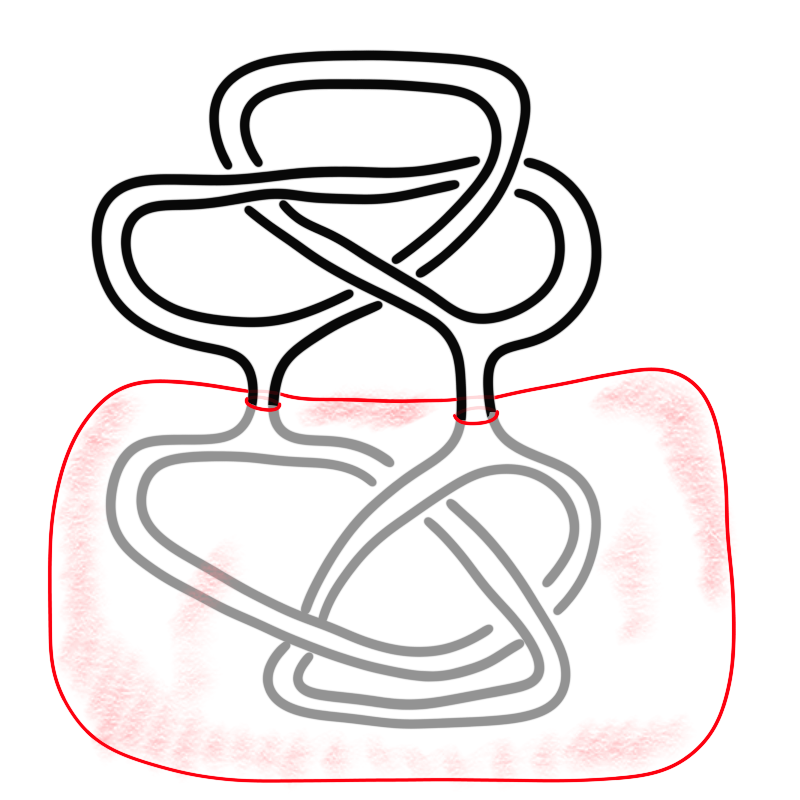}
\caption{}
\label{fig:Aone} 
\end{subfigure}
\begin{subfigure}{.3\linewidth}
	\centering 
	\includegraphics[scale=.1]{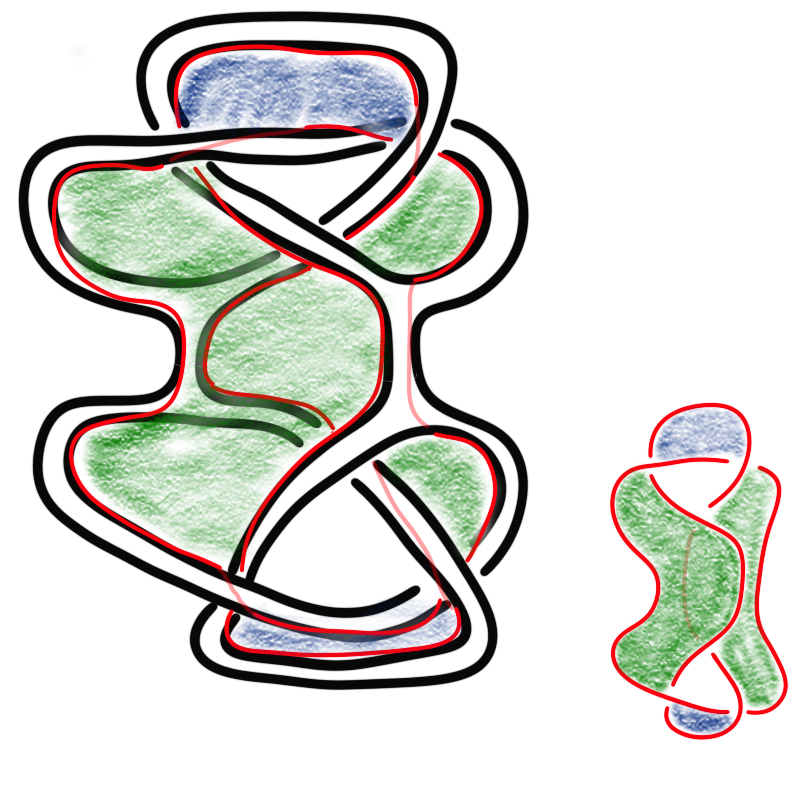} 
	\caption{}
	\label{fig:Atwo} 
\end{subfigure}
\begin{subfigure}{.3\linewidth}
	\centering 
	\includegraphics[scale=.13]{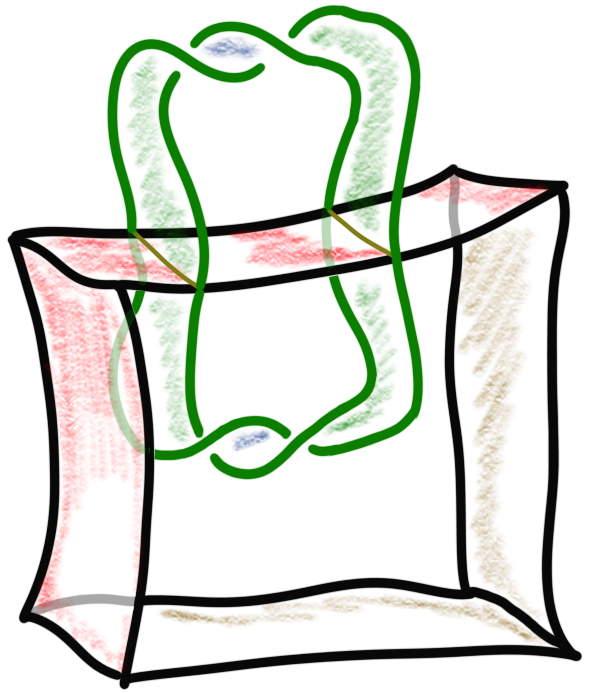} 
	\caption{}
	\label{fig:intersection} 
\end{subfigure}
\begin{subfigure}{.3\linewidth}
	\centering 
	\includegraphics[scale=.1]{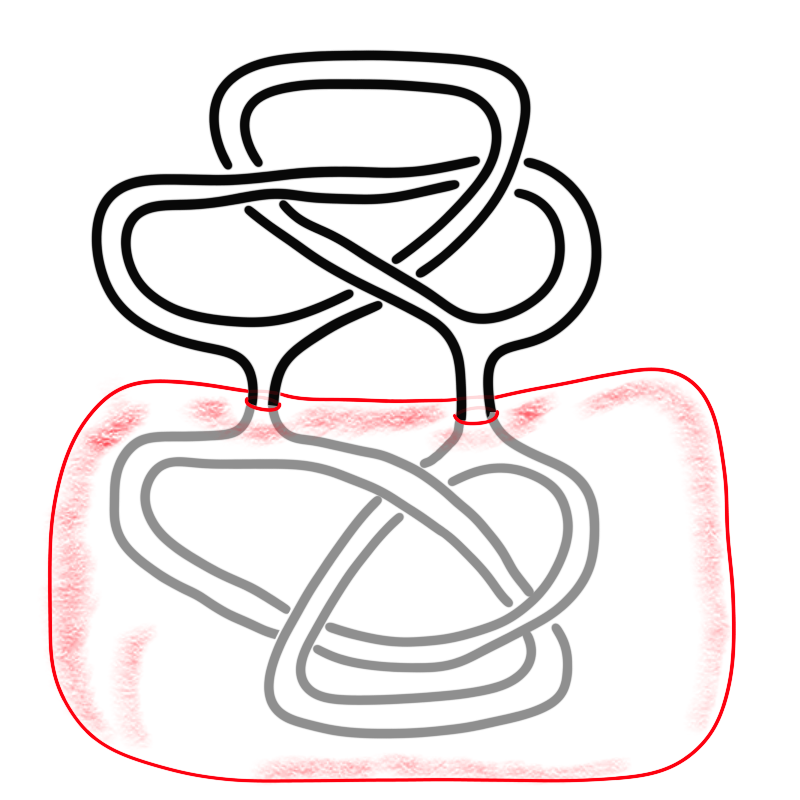} 
	\caption{}
	\label{fig:Aoneprime} 
\end{subfigure}
\begin{subfigure}{.3\linewidth}
	\centering 
	\includegraphics[scale=.1]{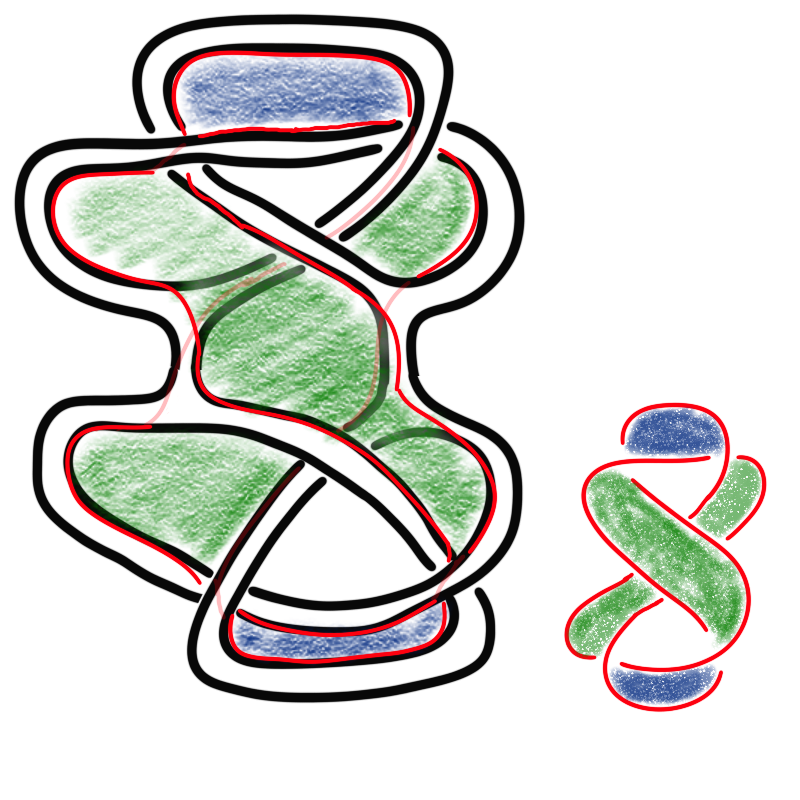} 
	\caption{}
	\label{fig:Atwoprime} 
\end{subfigure}
\begin{subfigure}{.24\linewidth}
	\centering 
	\includegraphics[scale=.13]{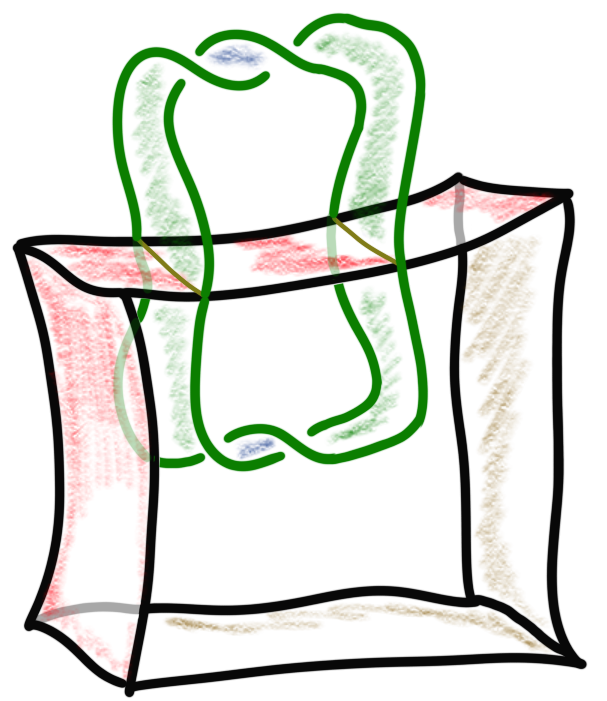} 
	\caption{}
	\label{fig:intersectionprime} 
\end{subfigure} 
\caption{Inequivalence of $6_{12}$ and $6_{13}$.}
\end{figure}

\begin{lemma}\label{lm:irre_hardpair}
$6_{12}$, $6_{13}$ are inequivalent, up to mirror image. 
\end{lemma}
\begin{proof}
Consider the spines of $6_{12},6_{13}$ in Table \ref{tab:hknots}, 
and observe, for each spine, that there is a $2$-sphere meeting the spine at two points and 
separating the spine in two parts, each of which has a constituent trefoil knot.
The $2$-sphere induces an essential annulus $A_1$ (resp.\ $A_1'$) in the exterior $E$ of $6_{12}$ (resp.\ the exterior $E'$ of $6_{13}$); see Figs.\ \ref{fig:Aone}, \ref{fig:Aoneprime}.

There is also a non-separating essential annulus 
$A_2$ (resp.\ $A_2'$) in $E$ (resp.\ $E'$); see Figs.\ \ref{fig:Atwo}, \ref{fig:Atwoprime}. The intersection $A_1\cap A_2$ (resp.\ $A_1'\cap A_2'$) consists of two essential arcs with opposite orientation; see Figs.\ \ref{fig:intersection}, \ref{fig:intersectionprime}. Hence a regular neighborhood $N$ of $A_1\cup A_2\subset E$ 
(resp.\ $N'$ of $A_1'\cup A_2'\subset E'$)
is an I-bundle over a twice-punctured Klein bottle.

In particular, $N$ (resp.\ $N'$) cuts off from 
$E$ (resp.\ $E'$) two solid tori $V_1,V_2$ (resp.\ $V_1'V_2'$). The frontier $F_i$ of $V_i$ (resp.\ $F_i'$ of $V_i'$) is precisely the cabling annulus of the constituent trefoil knot. Thus, the cores of $F_1,F_2$, with respect to $V_1,V_2$, have a slope of $\frac{2}{3}, -\frac{2}{3}$, while the cores of $F_1',F_2'$, with respect to $V_1',V_2'$, both have a slope of $\frac{2}{3}$. As a result $6_{12}$ is inequivalent to $6_{13}$ or its mirror image.  
\end{proof}  

Note that the above decomposition shows the exteriors of $6_{12},6_{13}$ have isomorphic fundamental groups, a reminiscence of the granny and square knots.

\subsection{A reducible hard pair}
\label{sec:typeone_hardpair}
\begin{figure}[ht]
\begin{subfigure}{.4\linewidth}
\centering
\includegraphics[scale=.6]{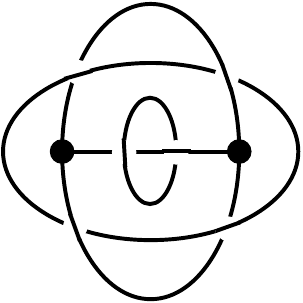} 
\caption{HL6$_{10}$.}
\label{fig:hl6_10}	
\end{subfigure}
\begin{subfigure}{.55\linewidth} 
\centering
\begin{overpic}[scale=.15,percent]{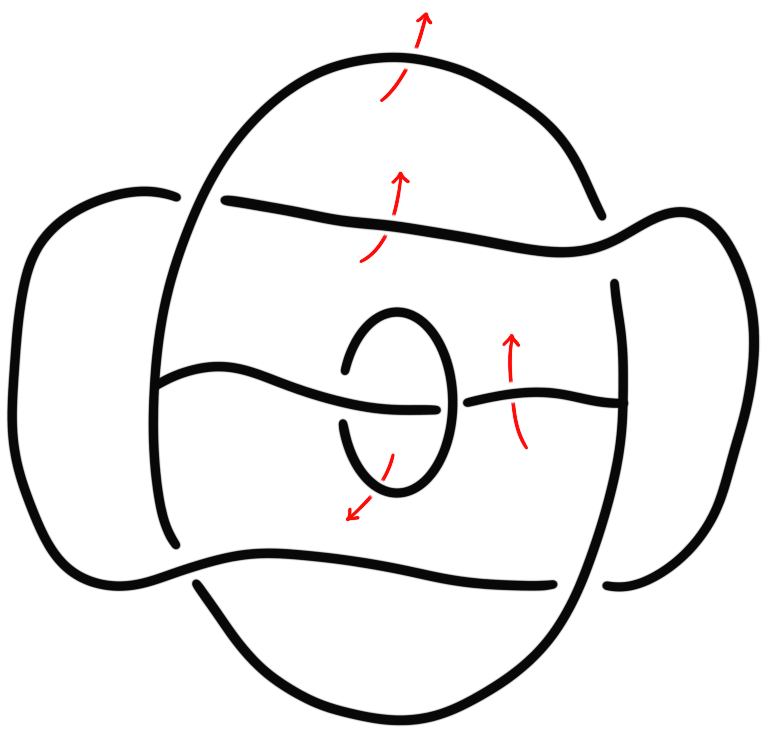}
\put (51,80){$a$}
\put (50,57){$b$}
\put (69,36){$c$}
\put (52,33){$d$}
\end{overpic}
\caption{Fundamental group of the exterior of HL$6_{10}$.}
\label{fig:wirtinger_hlsixten}
\end{subfigure}
%\caption{}
\end{figure}

Recall the handlebody-link HL$6_{10}$ is of type $[2,1]$; see Fig.\ \ref{fig:hl6_10}. The $1$-sum $V_1\in \mathrm{HL}6_{10}$ given by the genus two component is of type $[2,0,1]$, whereas the $1$-sums $V_2,V_3\in \mathrm{HL}6_{10}$ given by two genus one components of HL$6_{10}$ are of type $[1,2]$. The fundamental groups of their exteriors are all isomorphic, being the free product of $\mathbb Z$ and the fundamental group of the exterior of HL$6_{10}$.

\begin{lemma}\label{lm:re_hardpair}
$V_1,V_2,V_3$ are mutually inequivalent, up to mirror image.
\end{lemma}
\begin{proof}
$V_1$ can be distinguished from $V_2,V_3$ easily by its type. 
To see $V_2,V_3$ are inequivalent, up to mirror image, it suffices to show that there is no self-homeo\-morphism of $\sphere$ preserving HL$6_{10}$ and swapping the two genus one components. Denote by $W_2,W_3$ the genus one components of HL$6_{10}$, and let $l_2,l_3$ be their preferred longitudes. Attaching a $2$-handle along $l_i$ to the exterior of HL$6_{10}$ yields a $3$-manifold $M_i$, $i=2,3$. 

Any self-heomomorphim of $\sphere$ keeping HL$6_{10}$ and swapping the two genus one components induces a homeomorphism between $M_2$ and $M_3$. The fundamental groups of $M_2,M_3$ is the the fundamental group of the exterior of $\mathrm{HL}6_{10}$ with relators given by $l_2,l_3$ added. In terms of the generators in Fig.\ \ref{fig:wirtinger_hlsixten}, they are
\begin{align*}
\langle a,b,c, d & \mid\ dcd^{-1}c^{-1}=1, ba^{-1}b^{-1}c^{-1}a^{-1}b^{-1}acba=1,
 ba^{-1}b^{-1}c^{-1}a^{-1}=1\rangle,\\
\langle a,b,c, d & \mid\ dcd^{-1}c^{-1}=1, ba^{-1}b^{-1}c^{-1}a^{-1}b^{-1}acba=1,
dcd^{-1}=1\rangle.
\end{align*}
They have different $\texttt{ks}_{S_3}$-invariants, and hence are non-isomorphic, so $M_2,M_3$ are non-homeomorphic. Therefore, $V_2,V_3$ are inequivalent, up to mirror image. 
\end{proof} 

%%%%%%%%%%%%%%%%%%%%%%%%%%%%%%%%%%%%%%%%%%%%%%%%%%%%%%%%%%%%%%%%%%%%%%%%%%%%%

%%%%%%%%%%%%%%%%%%%%%%%%%%%%%%%%%%%%%%%%%%%%%%%%%%%%%%%%%%%%%%%%%%%%%%%%%%%%%
\section*{Appendix}\label{sec:appendix}
We include here an example of visual equivalence proofs via diagrams.
It shows diagrammatically the equivalence between the $12$ diagrams in
the bigtree final leaf corresponding to entry $5_1$, Table \ref{tab:hknots}.

Prior to this, the entries in this leaf were separated from all the remaining diagrams
in the bigtree by successively computing three increasingly powerful invariants,
namely the Kitano-Suzuki
invariant \texttt{ks}$_G$ with $G = A_4, A_5, S_6$, for example \texttt{ks}$_{A_5}$ does
not distinguish this entry from entry $6_8$, both in Table \ref{tab:hknots}.
For the record, the computed values of the three invariants are $210, 4988, 595165$ respectively.
The first two coincide when computed on $6_8$ whereas for $G = S_6$ the computation for entry
$6_8$ gives the value $601240$.

\begin{figure}
\includegraphics[width=0.8\textwidth]{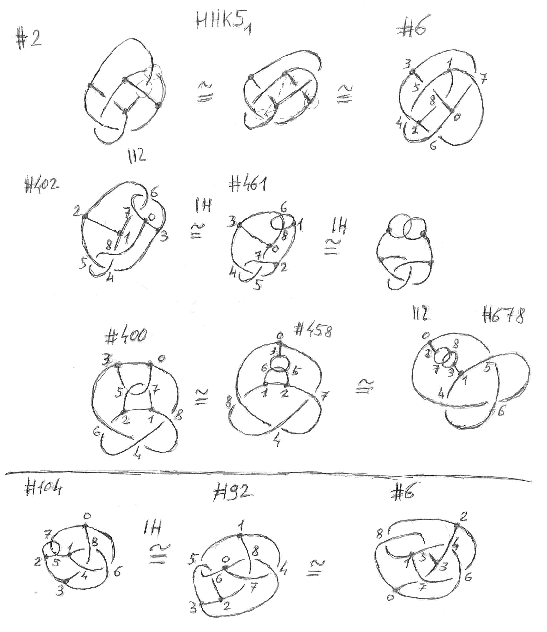}
\caption{Diagrammatic equivalences between entries in the same leaf as $5_1$ in Table \ref{tab:hknots}.}
\label{fig:equivs_HHK5_1}
\end{figure}

Figure \ref{fig:equivs_HHK5_1} shows mutual equivalence between all diagrams contained in the
considered leaf of the bigtree.
The number `\#$n$' attached to each diagram corresponds to the numbering assigned to the corresponding
plane graph by our software used to
produce all plane embeddings.
The overcrossing information is not indicated for \textit{loopings}, i.e.
those portions of the diagram with one trivalent node, two crosses and three
partial arcs, i.e. the result of a $3$-sum with entry $2_1$ in \cite[Table 1]{Mor:07b}.
There are two equivalent (implied) choices for those two overcrossings that do not immediately imply
nonminimality via a Reidemeister R2 move.

In three cases there are two different choices of the
overcrossing (at two crossings) for the same plane embedding.
This happens for entries \#92, \#400 and \#402.
The two different choices lead to equivalent diagrams as is apparent by applying an IH move followed by
a \texttt{loop-flip} (Fig. \ref{fig:autopairs}) and an IH move again.

Each displayed equivalence can be easily realized with short sequences of Reidemeister moves I to V plus
IH moves.
As an illustration the very first equivalence between two variants of diagram \#2 where two overcrossings
are reverted (corresponding to the \texttt{barrel-flip} situation in Figure \ref{fig:autopairs}) can
be realized with the sequence \texttt{inv(V)}, III, V.
The immediately following equivalence of the result with \#6 can be realized by first performing an \texttt{inv(IV)}
followed by a IV.

One particular equivalence, however, connecting the unnamed diagram on the second row with \#678, although rather
visually evident, entails a long sequence of extended Reidemeister moves R1-R5.

%%%%%%%%%%%%%%%%%%%%%%%%%%%%%%%%%%%%%%%%%%%%%%%%%%%%%%%%%%%%%%%%%%%%%%%%%%%%%
%
%%%%%%%%%%%%%%%%%%%%%%%%%%%%%%%%%%%%%%%%%%%%%%%%%%%%%%%%%%%%%%%%%%%%%%%%%%%%%%
%\section{Handlebody knots vs spatial graphs}
%\input{spatialgraphs}
%%%%%%%%%%%%%%%%%%%%%%%%%%%%%%%%%%%%%%%%%%%%%%%%%%%%%%%%%%%%%%%%%%%%%%%%%%%%%%
%
%%%%%%%%%%%%%%%%%%%%%%%%%%%%%%%%%%%%%%%%%%%%%%%%%%%%%%%%%%%%%%%%%%%%%%%%%%%%%%
%\section{Reducible handlebody knots}\label{sec:reducible}
%\input{reducible}
%%%%%%%%%%%%%%%%%%%%%%%%%%%%%%%%%%%%%%%%%%%%%%%%%%%%%%%%%%%%%%%%%%%%%%%%%%%%%%
%\input{tablered}
%%%%%%%%%%%%%%%%%%%%%%%%%%%%%%%%%%%%%%%%%%%%%%%%%%%%%%%%%%%%%%%%%%%%%%%%%%%%%%
%
%\section*{Acknowledgements}
%\textbf{Acknowledgements}.
%...
%The first author benefits from the support
%of the GNA\-MPA (Gruppo Nazionale per l'Analisi Matematica, la Probabilit\`a
%e le loro Applicazioni) of INdAM (Istituto Nazionale di Alta
%Matematica). The second author
%benefits from the support of the Swiss National Science Foundation Professorship grant PP00P2\_179110/1. 
%The fourth author is supported by
%National Center of Theoretical Sciences.

\section*{Acknowledgment}
%\textbf{Acknowledgements}.
The first author benefits from the support
of GNAMPA (INdAM). The fourth author gratefully acknowledges the support from NSTC (grant no. 115-2115-M-110 -002 -MY3).

\section*{Declaration of AI use}
After completing the manuscript we asked for an assessment of the content to GPT-6 Astra.
This allowed us to spot a couple of omissions.
Specifically, in Subsection \ref{sec:typetwospecial} we now also cover the case of a 2-sum performed at a solid torus
component of the $\fthree$ factor;
in Section \ref{sec:irreducibility} we now also show unsplittablity of all entries in Table \ref{tab:hlinks}.
These additions did not change any of the main results in the paper.

%%%%%%%%%%%%%%%%%%%%%%%%%%%%%%%%%%%%%%%%%%%%%%%%%%%%%%%%%%%%%%%%%%%%%%%%%%%%%

\end{document}